\documentclass[12pt]{article}
 
 \usepackage{amsmath,amssymb,amscd,amsthm,esint}
 
\usepackage{graphics,amsmath,amssymb,amsthm,mathrsfs}

\newtheorem{theorem}{Theorem}[section]
\newtheorem{prop}[theorem]{Proposition}
\newtheorem{lemma}[theorem]{Lemma}

\newtheorem{remark}[theorem]{Remark}

\def\xxint#1#2#3{{\setbox0=\hbox{$#1{#2#3}{\int}$}
  \vcenter{\hbox{$#2#3$}}\kern-.5\wd0}}

\def \rd {{\mathbb R}^d}
\def\rdd{{\mathbb R}^{d\times d}}

\def\e{\varepsilon}
\def\R{\mathbb{R}}
\def\cL{\mathcal{L}}

\newcommand{\average}{-\!\!\!\!\!\!\int}

\begin{document}

\title
{\bf Boundary Korn Inequality and \\ Neumann Problems  in Homogenization \\ of
Systems of Elasticity}

\author{Jun Geng\thanks{Supported in part by the NNSF of China (no. 11571152) and Fundamental Research Funds for the Central
Universities (LZUJBKY-2015-72).}\qquad 
Zhongwei Shen\thanks{Supported in part by NSF grant DMS-1600520.}
\qquad Liang Song \thanks{Supported in part by the NNSF of China (Nos 11471338 and 11622113) and Guangdong Natural Science Funds for Distinguished Young Scholar (No. 2016A030306040).}}


\date{ }

\maketitle

\begin{abstract}

This paper concerns with a family of elliptic systems of linear elasticity 
with rapidly oscillating periodic coefficients, arising in the theory of homogenization.
We establish uniform optimal  regularity estimates for solutions of Neumann problems
in a bounded Lipschitz domain with $L^2$ boundary data.
The proof relies on a boundary Korn inequality for solutions of systems of linear elasticity
and uses a large-scale Rellich estimate obtained in \cite{Shen-2016}.

\bigskip



\end{abstract}



\section{Introduction}
\setcounter{equation}{0}

This paper concerns with a family of elliptic systems of linear elasticity 
with rapidly oscillating periodic coefficients, arising in the theory of homogenization.
We establish uniform optimal  regularity estimates for solutions of Neumann problems
in a bounded Lipschitz domain with $L^2$ boundary data.
The proof relies on a boundary Korn inequality for solutions of systems of linear elasticity
and uses a large-scale Rellich estimate obtained in \cite{Shen-2016}.

More precisely, we consider a family of elasticity operators,
\begin{equation}\label{operator}
\cL_\e=-\text{\rm div} \big(A(x/\e)\nabla \big)=-\frac{\partial}{\partial x_i}
\left\{ a_{ij}^{\alpha\beta} \left(\frac{x}{\e}\right) \frac{\partial}{\partial x_j} \right\}, \quad \e>0,
\end{equation}
 where $A(y)=(a_{ij}^{\alpha\beta} (y))$
with $1\le i, j, \alpha, \beta\le d$.
Throughout this paper we will assume that the coefficient matrix (tensor) $A$ satisfies the
elasticity condition,
\begin{equation}\label{elasticity}
\aligned
 & a_{ij}^{\alpha\beta} (y) =a_{ji}^{\beta\alpha} (y) =a_{\alpha j}^{i\beta} (y),\\
 & \kappa_1 |\xi|^2 \le a_{ij}^{\alpha\beta}  (y) \xi_i^\alpha \xi_j^\beta \le \kappa_2 |\xi|^2
 \endaligned
 \end{equation}
for $y\in \R^d$ and for { symmetric} matrix $\xi=(\xi_i^\alpha)\in \R^{d\times d}$,
where $\kappa_1, \kappa_2>0$
(the summation convention is used throughout the paper).
We will also assume that $A(y)$ is 1-periodic, 
\begin{equation}\label{periodicity}
A(y+z)=A(y) \quad \text{ for  } y\in \R^d \text{ and } z\in \mathbb{Z}^d,
\end{equation}
and is H\"older continuous,
\begin{equation}\label{Holder}
\| A\|_{C^\sigma (\R^d)} \le M
\end{equation}
for some $\sigma \in (0,1)$ and $M>0$.

The following is the main result of the paper.

\begin{theorem}\label{main-theorem-2}
Assume that $A$ satisfies conditions (\ref{elasticity}), (\ref{periodicity}) and (\ref{Holder}).
Let $\Omega$ be a bounded Lipschitz domain in $\rd$.
Then for any $g\in L^2_{\cal{R}}(\partial\Omega)$, there exists a weak solution $u_\e$,
unique up to a rigid displacement, to the Neumann problem,
\begin{equation}\label{NP}
\left\{
\aligned
\mathcal{L}_\e (u_\e) & =0 & \ & \text{ in } \Omega,\\
 \frac{\partial u_\e}{\partial \nu_\e} & = g & \ & \text{ n.t. on } \partial \Omega,\\
 (\nabla u_\e)^* & \in L^2(\partial \Omega),
 \endaligned
 \right.
 \end{equation}
and the solution $u_\e$ satisfies the estimate
\begin{equation}\label{N-estimate}
\| (\nabla u_\e)^*\|_{L^2(\partial\Omega)} \le C\, \| g\|_{L^2(\partial\Omega)},
\end{equation}
where $C$ depends only on $d$, $\kappa_1$, $\kappa_2$, $(\sigma, M)$  and the Lipschitz 
character of $\Omega$.
\end{theorem}

We  introduce the notations used in Theorem \ref{main-theorem-2} and hereafter. 
We use 
$$
\frac{\partial u_\e}{\partial \nu_\e}=n\cdot A(x/\e)\nabla u_\e
$$
 to denote the conormal derivative  on $\partial\Omega$ of $u_\e$,
 associated with the operator $\cL_\e$, where $n$ is the outward unit normal to $\partial\Omega$.
 The boundary value in  (\ref{NP}) is taken a.e. in the sense of nontangential (n.t.) convergence.
 By $(w)^*$ we mean the nontangential maximal function of $w$, defined by
 \begin{equation}\label{nt}
 (w)^* (z) =\sup \Big\{ |w(x)|: \ x\in  \Omega \text{ and } |x-z|< C_0 \delta (x) \Big\}
 \end{equation}
for $z\in \partial\Omega$,
 where $\delta (x)=\text{dist}(x, \partial\Omega)$ and
 $C_0>1$ is a large constant depending on the Lipschitz character of $\Omega$.
Also, 
\begin{equation}\label{R}
\mathcal{R} =\big\{ \phi=Bx +b: B\in \R^{d\times d} \text{ is skew-symmetric  and } b\in \R^d \big\}
\end{equation}
is the space of rigid displacements and
\begin{equation}\label{L-R}
L^2_{\mathcal{R}} (\partial\Omega)
=\Big\{ g\in L^2(\partial\Omega; \R^d): \ \int_{\partial\Omega} g\cdot  \phi =0 \text{ for any } \phi\in \mathcal{R} \Big\}.
\end{equation}

Boundary value problems in Lipschitz domains with $L^p$ boundary data
have been studied extensively since late 1970's.
We refer the reader to the book \cite{Kenig-book} for references in this area up to mid-1990's
and to \cite{KKPT-2000, Verchota-2005, Shen-2006, Shen-2007, KS1, KS2, HKMP-2015, DPR-2016} 
and their references for more recent work.
For elliptic systems of linear elasticity, 
in  the case of a homogeneous isotropic body, where the coefficients 
are constants and given by
\begin{equation}\label{iso}
a_{ij}^{\alpha\beta}
=\mu \delta_{ij}\delta_{\alpha\beta} +\lambda \delta_{i\alpha} \delta_{j\beta}
+\mu \delta_{i\beta}\delta_{j\alpha}
\end{equation}
with Lam\'e constants $\lambda$ and $\mu$,
Theorem \ref{main-theorem-2} as well as the corresponding results for Dirichlet problem
was proved 
in \cite{DKV-1988} (also see \cite{Verchota-1984, fabes2, FKV-1988, Gao-1991, PV-1995}), by the method of layer potentials.
The general case of elliptic systems of elasticity with constant coefficients satisfying 
(\ref{elasticity}) was 
treated in \cite{Verchota-1986}.
Our estimate (\ref{N-estimate}) is new for variable coefficients, even in the (local) case $\e=1$.

For elliptic equations and systems with rapidly oscillating periodic coefficients,
the  Dirichlet and Neumann problems with $L^p$ boundary data
were studied in \cite{AL-1987-ho, AL-1987, AL-1991, KS1, KS2, KLS1}.
In particular, if $\Omega$ is Lipschitz,
the uniform estimate (\ref{N-estimate}) in Theorem \ref{main-theorem-2}
as well as the corresponding estimates
\begin{equation}\label{D-estimate-1}
\| (u_\e)^*\|_{L^2(\partial\Omega)} \le C\, \| f\|_{L^2(\partial\Omega)},
\end{equation}
\begin{equation}\label{R-estimate-1}
\| (\nabla u_\e)^*\|_{L^2(\partial\Omega)} \le C\, \| f\|_{W^{1,2}(\partial\Omega)},
\end{equation}
for solutions of Dirichlet problem:
$\mathcal{L}_\e (u_\e)=0$ in $\Omega$ and $u_\e =f$  on $\partial\Omega$,
was established in \cite{KS2},
where it is assumed that $A=(a_{ij}^{\alpha\beta})$,
with $1\le i, j\le d$ and $1\le \alpha, \beta\le m$,
is 1-periodic, symmetric, H\"older continuous, and satisfies the very
strong  ellipticity condition or the Legendre condition,
\begin{equation}\label{ellipticity}
\mu |\xi|^2 \le a_{ij}^{\alpha\beta} (y) \xi_i^\alpha\xi_j^\beta
\le \frac{1}{\mu} |\xi|^2 \quad 
\text{ for } y \in \R^d \text{ and  }
\xi=(\xi_i^\alpha) \in \R^{m\times d},
\end{equation}
where $\mu>0$.
We mention that for a scalar elliptic  equation $(m=1)$ in a Lipschitz domain, the estimate (\ref{D-estimate-1})
is due to B. Dahlberg (unpublished), while (\ref{R-estimate-1}) and (\ref{N-estimate}) were
proved in \cite{KS1}.
Also, if $\Omega$ is a bounded $C^{1, \alpha}$ domain and $1<p<\infty$,
the $L^p$ analogous of estimates (\ref{D-estimate-1}),
(\ref{N-estimate}) and (\ref{R-estimate-1}) may be found in \cite{AL-1987-ho, AL-1987, KLS1}
under the ellipticity condition (\ref{ellipticity}).

To describe the main difficulties  in the study of the $L^2$ Neumann problem (\ref{NP})
for systems of elasticity
and our approach to Theorem \ref{main-theorem-2} as well as the structure of this paper, we first note that
it is possible to rewrite the system $\text{\rm div}(A(x/\e)\nabla u_\e)=0$
as $\text{\rm div}(\widetilde{A}(x/\e)\nabla u_\e)=0$ in such a way that
the coefficient matrix $\widetilde{A}(y)$ is symmetric and satisfies (\ref{ellipticity}).
This allows us to use the interior Lipschitz estimates in \cite{AL-1987}
and the optimal estimates for Dirichlet problem in \cite{KS2}.
As a result,  estimates (\ref{D-estimate-1}) and (\ref{R-estimate-1}) hold for solutions of 
Dirichlet Problem for elliptic systems of elasticity (see Section 2).
We remark that the same technique was used in \cite{DKV-1988}
in the case of a homogeneous isotropic body. 

Although the rewriting of the system of elasticity $\mathcal{L}_\e (u_\e)=0$
does not change Dirichlet problem, it does change the Neumann problem as the conormal derivative 
$\partial u_\e/\partial\nu_\e$ depends on the coefficient matrix.
Nevertheless, it makes the method of layer potentials available to $\mathcal{L}_\e$,
since the estimates of fundamental solutions and layer potentials only involve the interior estimates.
As in the case of constant coefficients \cite{Verchota-1984, DKV-1988, FKV-1988,
Verchota-1986},
to use the method of layer potentials for $L^2$ Neumann  problems in Lipschitz domains,
the key step is  to establish the following Rellich estimate,
\begin{equation}\label{Rellich-estimate-N}
\int_{\partial\Omega} |\nabla u_\e|^2\, d\sigma 
\le C \int_{\partial\Omega} \Big| \frac{\partial u_\e}{\partial \nu_\e}\Big|^2 d\sigma
+\frac{C}{r_0} \int_\Omega |\nabla u_\e|^2\, dx,
\end{equation}
for suitable solutions of $\mathcal{L}_\e (u_\e)=0$ in $\Omega$,
where $r_0=\text{\rm diam}(\Omega)$. 
In fact, we will  show in Section 3 that 
for a given pair of $(\Omega, \e)$,
the estimate (\ref{Rellich-estimate-N}) is equivalent to (\ref{N-estimate}).

In comparison to elliptic systems with coefficients satisfying (\ref{ellipticity}) \cite{KS2},
one of the main difficulties in proving  (\ref{Rellich-estimate-N}) is caused by the fact that 
from the elasticity condition (\ref{elasticity}) one only obtains
\begin{equation}\label{coer}
\int_{\partial\Omega} A(x/\e)\nabla u_\e \cdot \nabla u_\e\, d\sigma
\ge c \int_{\partial\Omega} |\nabla u_\e + (\nabla u_\e)^T)|^2\, d\sigma,
\end{equation}
where $c>0$ and $(\nabla u_\e)^T$ denotes the transpose of the $d\times d$ matrix $\nabla u_\e$.
As a result, to control the full gradient $\nabla u_\e$ on $\partial\Omega$,
we would need  some Korn type inequality on the boundary.
We remark that the same issue already appears at the small scale, where $\e=1$ and
diam$(\Omega)\le 1$, even in the case of constant coefficients \cite{DKV-1988, Verchota-1986}.
 Also, the techniques developed in \cite{KS2} for treating the large-scale estimates
 fail due to the lack of (uniform) Korn inequalities on boundary layers.

Our proof of (\ref{Rellich-estimate-N}) is motivated by the work in \cite{DKV-1988, Verchota-1986, KS2}
and involves several new ideas. We divide the proof into five steps.
Note that the first four steps treat the case of small scales, in which 
the estimates are local and does not use the periodicity assumption.

Step 1. Let $\e=1$ and $\Omega$ be a Lipschitz domain with
$r_0=\text{\rm diam}(\Omega)\le (1/4)$.
We establish the boundary Korn inequality 
\begin{equation}\label{boundary-Korn-0}
\| \nabla u\|_{L^2(\partial\Omega)} 
\le C \, \| \nabla u +(\nabla u)^T\|_{L^2(\partial\Omega)}
+ C\,  r_0^{-1/2} \|\nabla u\|_{L^2(\Omega)}
\end{equation}
for solutions of $\mathcal{L}_1 (u)=0$ in $\Omega$, under the additional assumption:
\begin{equation}\label{additional-0}
|\nabla A(x)|\le M_0 \big[\delta (x)\big]^{\sigma -1} \quad \text{ and } \quad 
|\nabla^2  A(x)|\le M_0 \big[\delta (x)\big]^{\sigma -2}
\end{equation}
for any $x\in \Omega$, where $\sigma\in (0,1)$ and $M_0>0$.
To do this we rewrite the system $\mathcal{L}_1 (u)=0$ in $\Omega$ as
$$
\mu \Delta u^\alpha  +\frac{\partial}{\partial x_i} \left( b_{ij}^{\alpha\beta} \frac{\partial u^\beta}{\partial x_j} \right)
=0 \quad \text{ in } \Omega,
$$
where $B= (b_{ij}^{\alpha\beta})$ satisfies the elasticity condition (\ref{elasticity}), with (different) constants depending
on $\kappa_1$ and $\kappa_2$.
We then  factor the matrix $B(x) $ so that
$$
b_{ij}^{\alpha\beta} = q_i^{t\alpha}q _j^{t\beta},
$$
where $t$ is summed from $1$ to $m={d(d+1)}/{2}$.
Using the boundary Korn inequality (\ref{boundary-Korn-0}) for harmonic functions \cite{DKV-1988}
as well as Dahlberg's bilinear estimate \cite{Dahlberg-1986, Shen-2015}, the problem is reduced to the estimates of
the nontangential maximal function and the square function of 
$$
v=(v^t) =\left( q_j^{t\beta} \frac{\partial u^\beta}{\partial x_j} \right).
$$
To complete the step we observe that $v$ is a solution of an $m\times d$ elliptic system
$$
L (v)= F_0 + \text{\rm div}(F_1) \quad \text{ in } \Omega,
$$
with  coefficients satisfying  the Legendre condition (\ref{ellipticity}).
This allows us to use the estimates in \cite{KS2} to obtain the boundary Korn inequality 
 (\ref{boundary-Korn-0}).
 The details of the argument is given in Section 5, while some auxiliary estimates needed 
 for handling terms with $F_0$ and $F_1$ are given in Section 4.
 
 Step 2. 
Establish the Rellich estimate (\ref{Rellich-estimate-N}) for $\e=1$ and Lipschitz domain
 $\Omega$ with diam$(\Omega)\le (1/4)$,
 under the additional assumption (\ref{additional-0}).
 The proof uses a Rellich identity and
 (\ref{boundary-Korn-0}). See Section 6.
 
 Step 3. Let $\e=1$ and $\Omega$ be a Lipschitz domain with
diam$(\Omega)\le (1/4)$.
Prove (\ref{Rellich-estimate-N}) without the condition (\ref{additional-0}).
 This is done in Section 7, using an approximation scheme taken from  \cite{KS2}, with the help of 
 layer potentials.
 
 Step 4. Let $\mathcal{L}_\e (u_\e)=0$ in $\Omega$, where
  $0<\e< \text{\rm diam}(\Omega)<\infty$. It follows from Step 3 by some localization and rescaling techniques 
  that
  \begin{equation}\label{small-scale}
  \int_{\partial\Omega} |\nabla u_\e|^2\, d\sigma
  \le 
C \int_{\partial\Omega} \Big| \frac{\partial u_\e}{\partial \nu_\e}\Big|^2 d\sigma
+\frac{C}{\e } \int_{\Omega_{\e}}  |\nabla u_\e|^2\, dx,
\end{equation}
where $\Omega_r =\{ x\in \Omega: \text{dist} (x, \partial\Omega) < r \}$.

Step 5. To control the last term in (\ref{small-scale}), we use the following estimate 
\begin{equation}\label{large-scale}
\frac{1}{r} \int_{\Omega_r} |\nabla u_\e|^2\, dx
\le C  \int_{\partial\Omega} \Big| \frac{\partial u_\e}{\partial \nu_\e}\Big|^2 d\sigma
+\frac{C}{r_0 } \int_{\Omega}  |\nabla u_\e|^2\, dx,
\end{equation}
for any $\e\le r<r_0=\text{diam}(\Omega)$.
The estimate (\ref{large-scale}), which is due to homogenization of
$\mathcal{L}_\e$,
 was proved by the second author in \cite{Shen-2016}
for weak solutions of $\mathcal{L}_\e (u_\e)=0$ in $\Omega$
under the conditions (\ref{elasticity}) and (\ref{periodicity}), using a sharp convergence
rate in  $H^1(\Omega)$ for a two-scale expansion of $u_\e$.
Note that no smoothness condition  is  needed for (\ref{large-scale}).
We remark that estimate (\ref{large-scale}) may be regarded as a large-scale Rellich estimate for two reasons.
Firstly, by combining it with the small-scale estimate (\ref{small-scale}), we obtain the full Rellich estimate 
(\ref{Rellich-estimate-N}) and thus complete the proof of Theorem \ref{main-theorem-2} (see Section 8).
Secondly, if we are allowed to take the limit $r\to 0$ in (\ref{large-scale}), as in the case of constant coefficients,
we would recover the estimate (\ref{Rellich-estimate-N}).

We will use $C$ and $c$ to denote constants that depend at most on $d$, $\kappa_1$, $\kappa_2$, 
$(\sigma, M)$ in (\ref{Holder}), 
and the Lipschitz character of $\Omega$.
If a constant also depends on other parameters, it will be pointed out explicitly.
Finally, we shall use $\| g\|_2 $ to denote the norm of $g$ in $L^2(\partial\Omega)$.
The norm in $L^2(\Omega)$ will be denoted by $\| u\|_{L^2(\Omega)}$.


\section{Preliminaries}
\setcounter{equation}{0}

Let $A(y) =(a_{ij}^{\alpha\beta} (y))$, where $a_{ij}^{\alpha\beta} (y)$ satisfies the
elasticity condition (\ref{elasticity}).
Let
\begin{equation}\label{a-tilde}
\widetilde{a}_{ij}^{\alpha\beta} (y)=a_{ij}^{\alpha\beta} (y)
+\mu \delta_{i\alpha}\delta_{j\beta} -\mu \delta_{i\beta}\delta_{j\alpha}
\end{equation}
and
\begin{equation}\label{b}
b_{ij}^{\alpha\beta}  (y) =\widetilde{a}_{ij}^{\alpha\beta} (y)-\mu \delta_{ij}\delta_{\alpha\beta},
\end{equation}
where $\mu=\kappa_1/4>0$. Note that
\begin{equation}\label{symmetry-b}
\widetilde{a}_{ij}^{\alpha\beta}=\widetilde{a}_{ji}^{\beta\alpha}, 
\quad b_{ij}^{\alpha\beta}=b_{ji}^{\beta\alpha}=b_{\alpha j}^{i\beta}
\quad \text{ for any } 1\le i, j, \alpha, \beta\le d.
\end{equation}

\begin{prop}\label{prop-2.1}
Let $\widetilde{a}_{ij}^{\alpha\beta}$ be defined by (\ref{a-tilde}).
Then 
\begin{equation}\label{a-tilde-ellipticity}
\widetilde{a}_{ij}^{\alpha\beta} \xi_i^\alpha\xi_j^\beta
\ge \mu |\xi|^2 \qquad \text{ for  any } \xi =(\xi_i^\alpha)\in \mathbb{R}^{d\times d}.
\end{equation}
\end{prop}

\begin{proof}
Note that by (\ref{elasticity}),
\begin{equation}\label{elasticity-2}
a_{ij}^{\alpha\beta} \xi_i^\alpha\xi_j^\beta\ge \frac{\kappa_1}{4} |\xi +\xi^T|^2 \qquad 
\text{ for any } \xi=(\xi_i^\alpha)\in \rdd,
\end{equation}
where $\xi^T$ denotes the transpose of $\xi$.
It follows that
$$
\aligned
b_{ij}^{\alpha\beta} \xi_i^\alpha\xi_j^\beta
&= a_{ij}^{\alpha\beta} \xi_i^\alpha\xi_j^\beta -\mu |\xi|^2 +\mu \xi_i^i \xi_j^j -\mu \xi_i^j \xi_j^i\\
&\ge \frac{\kappa_1}{4} |\xi +\xi^T|^2
-\mu |\xi|^2 -\mu \xi_i^j \xi_j^i\\
&=\frac{\kappa_1}{2} \big( |\xi|^2 + \xi_i^\alpha \xi_\alpha^i \big)-\mu |\xi|^2 -\mu \xi_i^j \xi_j^i\\
&= \frac12 \Big(\frac{\kappa_1}{2}-\mu\Big) |\xi +\xi^T|^2.
\endaligned
$$
Since $\mu=\kappa_1/4$, we obtain 
\begin{equation}\label{b-ellipticity}
b_{ij}^{\alpha\beta} \xi_i^\alpha\xi_j^\beta \ge \frac{\kappa_1}{8} |  \xi +\xi^T|^2
\qquad \text{ for any } \xi= (\xi_i^\alpha) \in \rdd,
\end{equation}
and 
$$
\widetilde{a}_{ij}^{\alpha\beta}\xi_i^\alpha\xi_j^\beta
=\mu|\xi|^2 + b_{ij}^{\alpha\beta}\xi_i^\alpha\xi_j^\beta
\ge \mu |\xi|^2
$$
for any $\xi =(\xi_i^\alpha)\in \rd$.
\end{proof}

\begin{prop}\label{prop-2.2}
 Let
$\widetilde{A} (y)= (\widetilde{a}_{ij}^{\alpha\beta} (y))$ and
$u_\e\in H_{loc}^1(\Omega; \rd)$. Then 
$$
\text{\rm div} \big( A(x/\e)\nabla u_\e\big)=F \quad \text{ in }
\Omega \quad \text{  if and only if } \quad 
\text{\rm div} \big( \widetilde{A}(x/\e)\nabla u_\e\big)=F
\quad \text{ in }\Omega,
$$
where $F\in \big(C_0^\infty(\Omega; \rd)\big)^\prime$ is a distribution.
\end{prop}

\begin{proof}
Let $u=(u^\alpha) \in C^\infty (\Omega; \rd)$ and $ \varphi=(\varphi^\alpha)\in C_0^\infty(\Omega; \rd)$.
It follows from integration by parts that
\begin{equation}\label{2.2-1}
\aligned
\int_\Omega \Big\{ \delta_{i\alpha}\delta_{j\beta} -\delta_{i\beta} \delta_{j\alpha} \Big\}
\frac{\partial u^\beta}{\partial x_j} \cdot \frac{\partial \varphi^\alpha}{\partial x_i} 
&= \int_\Omega \text{\rm div} (u) \cdot \text{\rm div} (\varphi)
-\int_\Omega \frac{\partial u^\beta}{\partial x_\alpha} \cdot \frac{\partial \varphi^\alpha}{\partial x_\beta}\\
&= 0.
\endaligned
\end{equation}
By a density argument   (\ref{2.2-1}) continues to hold for any
$u\in H_{loc}^1(\Omega; \rd)$ and $\varphi\in C^\infty_0(\Omega; \rd)$.
Hence,
$$
\int_\Omega A(x/\e)\nabla u_\e \cdot \nabla \varphi\, dx
-\int_\Omega \widetilde{A}(x/\e)\nabla u_\e \cdot \nabla \varphi\, dx =0
$$
for any $u_\e \in H_{loc}^1(\Omega; \rd)$ and $\varphi\in C^\infty_0(\Omega; \rd)$.
\end{proof}

Proposition \ref{prop-2.2} reduces interior estimates for the operator $\mathcal{L}_\e$
to those for $\widetilde{\mathcal{L}}_\e$, where 
\begin{equation}\label{tilde-L}
\widetilde{\mathcal{L}}_\e =-\text{div} (\widetilde{A}(x/\e)\nabla ).
\end{equation}
Note that
by Proposition \ref{prop-2.1},
the coefficient matrix $\widetilde{A}$ satisfies the very strong ellipticity condition (\ref{ellipticity}).
It follows directly from \cite{AL-1987} that if $\mathcal{L}_\e (u_\e)=0$
in $B(x_0, r)$, then
\begin{equation}\label{interior-Lip}
|\nabla u_\e (x_0)|\le 
C \left(\average_{B(x_0, r)} |\nabla u_\e|^2 \right)^{1/2},
\end{equation}
where $C$ depends only on $\kappa_1$, $\kappa_2$, and $(\sigma, M)$.
Proposition \ref{prop-2.2} also shows that Dirichlet problem for $\mathcal{L}_\e$ is the same as
that for $\widetilde{\mathcal{L}}_\e$.
Since $\widetilde{A}$ is very strongly  elliptic, symmetric, and H\"older continuous,
the results in  \cite{KS2} for Dirichlet problems gives the following theorem.

\begin{theorem}\label{D-theorem}
Suppose that $A$ satisfies conditions (\ref{elasticity}), (\ref{periodicity}) and 
(\ref{Holder}).
Let $\Omega$ be a bounded Lipschitz domain in $\R^d$.
Then for any $f\in L^2(\partial\Omega; \R^d)$, there exists a unique solution $u_\e$
to Dirichlet problem, 
\begin{equation}\label{DP}
\left\{
\aligned
\mathcal{L}_\e (u_\e) & =0 & \ & \text{ in } \Omega,\\
 u_\e & = f & \ & \text{ n.t. on } \partial \Omega,\\
 (u_\e)^* & \in L^2(\partial \Omega),
 \endaligned
 \right.
 \end{equation}
and $u_\e$ satisfies the estimate (\ref{D-estimate-1}).
Furthermore, if $f \in W^{1,2}(\partial\Omega; \R^d)$, then the solution 
satisfies (\ref{R-estimate-1}).
The constants $C$ in (\ref{D-estimate-1}) and (\ref{R-estimate-1})
 depend only on $\kappa_1$, $\kappa_2$, $(\sigma, M)$, and the Lipschitz 
character of $\Omega$.
\end{theorem}

We end this section with a Rellich estimate for $\mathcal{L}_\e$ in $\Omega_-=\rd \setminus \overline{\Omega}$.
We will use $\nabla_{\tan} u$ to denote the tangential gradient of $u$ on $\partial\Omega$.

\begin{theorem}\label{Rellich-theorem-1}
Assume  that $A$ and $\Omega$ satisfy the same conditions as in Theorem \ref{D-theorem}.
Suppose that $\mathcal{L}_\e (u_\e) =0 $ in $\Omega_-$,
$(\nabla u_\e)^*\in L^2(\partial\Omega)$, and $\nabla u_\e $ exists n.t. on $\partial\Omega$.
Then
\begin{equation}\label{Rellich-estimate-1}
\int_{\partial\Omega} |\nabla u_\e|^2\,  d\sigma
\le C \int_{\partial\Omega} |\nabla_{\tan} u_\e|^2\, d\sigma
+\frac{C}{r_0} \int_{\Omega_-}  |\nabla u_\e|^2\, dx,
\end{equation}
where $r_0=\text{\rm diam}(\Omega)$ and $C$ depends only on
$\kappa_1$, $\kappa_2$, $(\sigma, M)$, and the Lipschitz character of $\Omega$.
\end{theorem}

\begin{proof}
Fix $z\in \partial\Omega$. Since $\widetilde{\mathcal{L}}_\e (u_\e)=0$ in $B(z, r_0)\cap \Omega_-$,
it follows from \cite{KS2} that
$$
\int_{B(z, cr_0)\cap \partial\Omega} |\nabla u_\e|^2\, d\sigma
\le C \int_{B(z, 2c r_0)\cap \partial\Omega} |\nabla_{\tan} u_\e|^2 \, d\sigma
+\frac{C}{r_0} \int_{B(z, 2c r_0)\cap \Omega_-} |\nabla u_\e|^2\, dx.
$$
Estimate (\ref{Rellich-estimate-1}) follows from this by 
covering $\partial \Omega$ with a finite number of balls
$\{ B(z_j, cr_0)\}$ centered on $\partial\Omega$.
\end{proof}



\section{Method of layer potentials}
\setcounter{equation}{0}

In \cite{KS2} the $L^2$ Dirichlet and Neumann problems in Lipschitz domains 
for elliptic systems with rapidly oscillating periodic coefficients
satisfying (\ref{ellipticity}) are solved by the 
method of layer potentials.
To solve the $L^2$ Neumann problem for the elasticity operator $\mathcal{L}_\e$,
we will also use the method of layer potentials.
Let $\Gamma_\e (x,y)=(\Gamma_\e^{\alpha\beta} (x,y))$ denote the $d\times d$ matrix
of fundamental solutions for $\mathcal{L}_\e$ (and $\widetilde{\mathcal{L}}_\e$).
It follows from the interior Lipschitz estimate (\ref{interior-Lip}) that
\begin{equation}\label{f-solution}
\aligned
|\Gamma_\e (x, y)| &\le C\,  |x-y|^{2-d},\\
|\nabla_ x\Gamma (x,y)| +|\nabla_y \Gamma_\e (x, y)| & \le C\,  |x-y|^{1-d},\\
|\nabla_x \nabla_y \Gamma_\e (x, y)| & \le C\,  |x-y|^{-d},
\endaligned
\end{equation}
for any $x, y\in  \rd$ and $x\neq y$ (some modifications are needed for the first estimate in the case $d=2$) .
Let $u_\e =\mathcal{S}_\e (f)$, where $f\in L^2(\partial\Omega; \rd)$ and 
\begin{equation}\label{single-layer}
\mathcal{S}_\e (f) (x)
=\int_{\partial\Omega} \Gamma_\e (x, y) f(y)\, d\sigma (y)
\end{equation}
denotes the single layer potential with density $f$.
Then
$\mathcal{L}_\e (u_\e)=0$ in $\rd\setminus \partial\Omega$
and $|u_\e (x)|= O(|x|^{2-d})$, $ |\nabla u_\e (x)|= O(|x|^{1-d})$
as $|x|\to \infty$.

We will use $(\nabla u_\e)_\pm$ to denote the nontangential limits of $\nabla u_\e$
taken from $\Omega_\pm$ respectively, where $\Omega_+=\Omega$
and $\Omega_- =\rd\setminus \overline{\Omega}$.

\begin{theorem}\label{layer-potential-theorem}
Suppose that $A$ satisfies (\ref{elasticity}), (\ref{periodicity}) and (\ref{Holder}).
Let $\Omega$ be a bounded Lipschitz domain in $\rd$.
Let $u_\e =\mathcal{S}_\e (f)$, where $f\in L^2(\partial\Omega; \rd)$. Then
\begin{equation}\label{layer-1}
\|(\nabla u_\e)^*\|_2 \le C\, \| f\|_2,
\end{equation}
where $C$ depends only on $\kappa_1$, $\kappa_2$, $(\sigma, M)$, and the Lipschitz character of $\Omega$.
Also, $(\nabla u_\e)_{\pm}$ exists n.t. on $\partial\Omega$, 
\begin{equation}\label{nt-limit}
(\nabla_{\tan} u_\e )_+ = (\nabla_{\tan}  u_{\e})_- \quad \text{ n.t.  on } \partial\Omega,
\end{equation}
and
\begin{equation}\label{co-normal-limit}
\left(\frac{\partial u_\e}{\partial \nu_\e}\right)_\pm
=\Big( \pm \frac12  I + \mathcal{K}_\e \Big) f \quad \text{ n.t. on } \partial\Omega,
\end{equation}
where $\mathcal{K}_\e$ is a bounded linear operator on $L^2(\partial\Omega; \rd)$.
\end{theorem}

\begin{proof} Estimates (\ref{layer-1}) and (\ref{nt-limit}) follow directly from  \cite[Section 4]{KS2}, as $\mathcal{L}_\e$
and $\widetilde{\mathcal{L}}_\e$ share the same fundamental solutions in $\rd$.
To see (\ref{co-normal-limit}), one uses  \cite[Theorem 4.4]{KS2} and  also the fact that
$$
a_{ij}^{\alpha\beta} (y) \eta_i \eta_j
=\widetilde{a}_{ij}^{\alpha\beta} (y) \eta_i\eta_j
$$ 
for any $ \eta =(\eta_i)\in \rd$.
\end{proof}

Note that $A(x/\e)\nabla v=0$ in $\rd$ for any $v\in \mathcal{R}$.
 Thus, if  $u_\e =\mathcal{S}_\e (f)$ and $v\in \mathcal{R}$, then
$$
\int_{\partial\Omega} \frac{\partial u_\e}{\partial\nu_\e} \cdot v\, d\sigma
=\int_\Omega A(x/\e)\nabla u_\e \cdot \nabla v\, dx =0.
$$
It follows from this and (\ref{layer-1}) that the operator 
\begin{equation}\label{L-2-R-1}
(1/2) I +\mathcal{K}_\e : L^2_\mathcal{R} (\partial\Omega)
\to L^2_\mathcal{R} (\partial\Omega) \quad \text{ is bounded uniformly in $\e>0$.}
\end{equation}
One of the main goals of this paper is to prove the following theorem,
which will allow us to solve the $L^2$ Neumann problem (\ref{NP}) for $\mathcal{L}_\e$
with optimal estimates.

\begin{theorem}\label{inverse-theorem}
Suppose that $A$ satisfies (\ref{elasticity}), (\ref{periodicity}) and (\ref{Holder}).
Let $\Omega$ be a bounded Lipschitz domain in $\rd$.
Then $(1/2)I +\mathcal{K}_\e: L^2_\mathcal{R} (\partial\Omega) \to L^2_\mathcal{R} (\partial\Omega)$
is invertible and
\begin{equation}\label{inverse-bound}
\| f\|_2 \le C \, \| \big((1/2)I +\mathcal{K}_\e \big) f\|_2 \quad \text{ for any } f\in L^2_\mathcal{R} (\partial\Omega),
\end{equation} 
where $C$ depends only on $\kappa_1$, $\kappa_2$, $(\sigma, M)$, and
the Lipschitz character of $\Omega$.
\end{theorem}

We now give the proof of Theorem \ref{main-theorem-2}, assuming Theorem \ref{inverse-theorem}.

\begin{proof}[\bf Theorem \ref{inverse-theorem} $\Longrightarrow$ Theorem \ref{main-theorem-2}]
The uniqueness follows directly from Green's identity,
\begin{equation}\label{Green-00}
\int_\Omega 
A(x/\e)\nabla u_\e \cdot \nabla u_\e  \, dx
=\int_{\partial\Omega}  u_\e \cdot \frac{\partial u_\e}{\partial \nu_\e}\, d\sigma.
\end{equation}
To show the existence as well as the estimate (\ref{N-estimate}),
for any $g\in L^2_{\mathcal{R}} (\partial\Omega)$, we choose $f\in L^2_{\mathcal{R}} (\partial\Omega)$
such that
$g=((1/2) I + \mathcal{K}_\e) f$ on $\partial\Omega$.
Let $u_\e =\mathcal{S}_\e (f)$. Then $u_\e$ is a solution to (\ref{NP}) and
$$
\|(\nabla u_\e)^*\|_2 \le C\,  \| f\|_2 \le C\,  \| g\|_2,
$$
where the last inequality follows from Theorem \ref{inverse-theorem}.
\end{proof}

\begin{remark}\label{remark-3.0}
{\rm Let $A$ be a matrix satisfying (\ref{elasticity}), (\ref{periodicity}) and (\ref{Holder}).
Let $\Omega$ be a bounded Lipschitz domain in $\rd$.
Suppose that the conclusions of Theorem \ref{main-theorem-2} hold for $A$ and $\Omega$.
Let $u_\e$ be a weak solution of $\mathcal{L}_\e (u_\e)=0$ in $\Omega$
such that $(\nabla u_\e)^* \in L^2(\partial\Omega)$ and
$\nabla u_\e$ exists n.t. on $\partial\Omega$.
Let $v_\e$ be the weak solution of $\mathcal{L}_\e (v_\e)=0$ in $\Omega$, given by Theorem \ref{main-theorem-2},
with Neumann data
$\partial v_\e/{\partial \nu_\e} =\partial u_\e/\partial \nu_\e$ on $\partial\Omega$.
Then 
\begin{equation}\label{remark-3.0-0}
\|(\nabla v_\e)^*\|_2 \le C\,  \Big\|\frac{\partial u_\e}{\partial \nu_\e}\Big \|_2.
\end{equation}
Since  $u_\e -v_\e =\phi \in \mathcal{R}$ in $\Omega$ and $\nabla \phi$ is constant, it follows that
\begin{equation}\label{remark-3.0-1}
\aligned
\int_{\partial\Omega} |\nabla u_\e|^2\,  d\sigma
&\le 2 \int_{\partial\Omega} |\nabla v_\e |^2 \, d \sigma 
+2 \int_{\partial\Omega} |\nabla \phi|^2 \, d\sigma  \\
& \le C \int_{\partial\Omega} \Big|\frac{\partial u_\e}{\partial \nu_\e}\Big|^2 \, d\sigma
+ \frac{C}{r_0^{d+1}} \Big| \int_\Omega \nabla \phi\, dx \Big|^2\\
&\le  C \int_{\partial\Omega} \Big|\frac{\partial u_\e}{\partial \nu_\e}\Big|^2 \, d\sigma
+ \frac{C}{r_0^{d+1}} \Big| \int_\Omega \nabla v_\e\, dx  \Big|^2
+\frac{C}{r_0^{d+1}} \Big| \int_\Omega \nabla u_\e  \, dx \Big|^2,
\endaligned
\end{equation}
where $r_0 =\text{diam}(\Omega)$. This,  together with (\ref{remark-3.0-0}),
gives
\begin{equation}\label{Rellich-estimate-0}
\int_{\partial\Omega} |\nabla u_\e|^2\,  d\sigma
\le C \int_{\partial\Omega} \left|\frac{\partial u_\e}{\partial \nu_\e}\right|^2 d\sigma
+\frac{C}{r_0^{d+1}} \Big| \int_\Omega \nabla u_\e  \, dx \Big|^2.
\end{equation}
In particular, for a given triple  $(\e, A, \Omega)$,
the invertibility  of $(1/2)I +\mathcal{K}_\e$ on $L^2_{\mathcal{R}} (\partial\Omega)$
and estimate (\ref{inverse-bound}) imply the Rellich estimate  (\ref{Rellich-estimate-0}).
}
\end{remark}

Next we will show that the Rellich estimate (\ref{Rellich-estimate-0})
implies (\ref{inverse-bound}).
To do this we  need some estimates for volume integrals.

\begin{lemma}\label{lemma-3.10}
Let $u_\e$ be a weak solution of $\mathcal{L}_\e (u_\e)=0$ in $\Omega$.
Suppose that $(\nabla u_\e)^*\in L^2(\partial\Omega)$ and $\nabla u_\e$ exists n.t. on $\partial\Omega$.
Then
\begin{equation}\label{L-0}
\| \nabla u_\e + (\nabla u_\e)^T\|_{L^2(\Omega)}
\le C\,  r_0^{1/2}\,  \Big\| \frac{\partial u_\e}{\partial \nu_\e} \Big\|_2,
\end{equation}
where $r_0=\text{\rm diam}(\Omega)$
and $\Omega$ depends only on $\kappa_1$, $\kappa_2$, and the Lipschitz character of $\Omega$.
\end{lemma}

\begin{proof}
By rescaling we may assume that $r_0=1$.
It follows by Green's identity that
\begin{equation}\label{Green-0}
\int_\Omega 
A(x/\e)\nabla u_\e 
\cdot \nabla u_\e  \, dx
=\int_{\partial\Omega} (u_\e  -v) \cdot \frac{\partial u_\e}{\partial \nu_\e}\, d\sigma,
\end{equation}
for any $v\in \mathcal{R}$. Hence, by (\ref{elasticity-2}) and the Cauchy inequality, 
\begin{equation}\label{L-0-1}
\frac{\kappa_1}{4} \int_\Omega |\nabla u_\e + (\nabla u_\e)^T|^2\, dx
 \le \| u_\e -v\|_2 \Big\| \frac{\partial u_\e}{\partial \nu_\e}\Big\|_2.
 \end{equation}
Using the well known trace inequality,
\begin{equation}\label{trace}
\int_{\partial \Omega} |u|^2 \, d\sigma 
\le C \int_\Omega |\nabla u|^2 \, dx +C \int_\Omega |u|^2 \, dx 
\end{equation}
for $u\in H^1(\Omega)$, and the second Korn inequality \cite{OSY-1992},
\begin{equation}\label{second-Korn}
\int_\Omega |\nabla u|^2\, dx +\int_\Omega |u|^2\, dx
\le C \int_\Omega |\nabla u + (\nabla u)^T|^2\, dx,
\end{equation}
which holds for any $u\in H^1(\Omega; \rd)$ such that $u\perp \mathcal{R}$ in $L^2(\Omega; \rd)$,
we see that
\begin{equation}\label{L-0-2}
\| u_\e -v\|_2  \le C\,  \| \nabla u_\e +(\nabla u_\e)^T\|_{L^2(\Omega)},
\end{equation}
if $v\in \mathcal{R}$ is chosen so that $u_\e -v \perp \mathcal{R}$ in $L^2(\Omega; \rd)$.
This, together with (\ref{L-0-1}), gives (\ref{L-0}).
\end{proof}

\begin{lemma}\label{lemma-3.20}
Let $u_\e =\mathcal{S}_\e (f)$, where $f\in L^2_\mathcal{R} (\partial\Omega) $.
Then
\begin{equation}\label{L-1}
\int_{\Omega_+ } |\nabla u_\e|^2\, dx +\int_{\Omega_- } |\nabla u_\e|^2\, dx
\le C\,  r_0  \| f\|_2 \Big\| \left(\frac{\partial u_\e}{\partial \nu_\e}\right)_+ \Big\|_2,
\end{equation}
where $r_0=\text{\rm diam}(\Omega)$
and $C$ depends only on $\kappa_1$, $\kappa_2$, and the Lipschitz character of $\Omega$.
\end{lemma}

\begin{proof}
By rescaling we may assume that $r_0=1$.
Let $ {\partial u_\e}/{\partial \widetilde{\nu}_\e}$ denote the conormal derivative on $\partial\Omega$ of
$u_\e$ with respect to the operator $\widetilde{\mathcal{L}}_\e$, given by (\ref{tilde-L}).
Since $\widetilde{\mathcal{L}}_\e (u_\e) =0$ in $\rd \setminus \partial\Omega$,
it follows from integration by parts that
\begin{equation}\label{Green}
\int_{\Omega_\pm}
\widetilde{A}(x/\e)\nabla u_\e \cdot \nabla u_\e \, dx
=\pm \int _{\partial\Omega}  u_\e \left( \frac{\partial u_\e}{\partial \widetilde{\nu}_\e}\right)_\pm \, d\sigma.
\end{equation}
This, together with the jump relation
\begin{equation}\label{jump-relation}
f=\left( \frac{\partial u_\e}{\partial \widetilde{\nu}_\e}\right)_+ 
-\left( \frac{\partial u_\e}{\partial \widetilde{\nu}_\e}\right)_- \quad \text{ on } \partial\Omega,
\end{equation}
gives
\begin{equation}\label{L-10}
\int_{\Omega_+}
\widetilde{A} (x/\e)\nabla u_\e \cdot \nabla u_\e \, dx
+\int_{\Omega_-}
\widetilde{A}(x/\e) \nabla u_\e \cdot \nabla u_\e \, dx
=\int_{\partial\Omega} u_\e \cdot f \, d\sigma.
\end{equation}
Hence, by the very strong  ellipticity of $\widetilde{A}$,
\begin{equation}\label{L-11}
\aligned
\int_{\Omega_+}  |\nabla u_\e|^2\, dx
+\int_{\Omega_-}  |\nabla u_\e|^2\, dx
 & \le C\,  \Big| \int_{\partial\Omega} u_\e \cdot f \, d\sigma \Big|\\
 & \le C \, \| u_\e -v\|_2 \| f\|_2,
 \endaligned
\end{equation}
for any $v\in \mathcal{R}$, where we have used the fact $f\in L^2_{\mathcal{R}} (\partial\Omega)$.
As in the proof of the last proposition, we may choose $v\in \mathcal{R}$ so that
$$
\| u_\e -v\|_2 \le C\,  \|\nabla u_\e + (\nabla u_\e)^T\|_{L^2(\Omega)}
\le C\,  \Big\|\left( \frac{\partial u_\e}{\partial \nu_\e}\right)_+ \Big\|_2,
$$
where we have used (\ref{L-0}) for the last step.
This, together with (\ref{L-11}), yields (\ref{L-1}).
\end{proof}

\begin{remark}\label{remark-3.20}
{\rm 
Let $A$ be a matrix satisfying (\ref{elasticity}), (\ref{periodicity}) and (\ref{Holder}).
Let $\Omega$ be a bounded Lipschitz domain in $\rd$.
Suppose that for any weak solution $u_\e$ of $\mathcal{L}_\e (u_\e)=0$ in $\Omega$ with
the properties that $(\nabla u_\e)^*\in L^2(\partial\Omega)$ and $\nabla u_\e$ exists
n.t. on $\partial\Omega$, the Rellich estimate  (\ref{Rellich-estimate-0}) holds.
We claim that this implies the estimate (\ref{inverse-bound}) for the operator 
$(1/2)I +\mathcal{K}_\e$. Indeed, let $u_\e =\mathcal{S}_\e (f)$ for some $f\in L^2_{\mathcal{R}} (\partial\Omega)$.
It follows from (\ref{Rellich-estimate-1}) and (\ref{Rellich-estimate-0}) as well as the jump
relation 
\begin{equation}\label{jump}
f =\left(\frac{\partial u_\e}{\partial \nu_\e}\right)_+
-\left(\frac{\partial u_\e}{\partial \nu_\e}\right)_- \quad \text{ on } \partial\Omega,
\end{equation}
that 
$$
\aligned
\| f\|_2  & \le \Big\| \left(\frac{\partial u_\e}{\partial \nu_\e}\right)_+\Big\|_2 
+ \Big\| \left(\frac{\partial u_\e}{\partial \nu_\e}\right)_-\Big\|_2 \\
& \le  \Big \| \left(\frac{\partial u_\e}{\partial \nu_\e}\right)_+\Big\|_2 
+ C \|\nabla_{\tan} u_\e \|_2 + C r_0^{-1/2} \| \nabla u_\e\|_{L^2(\Omega_-)} \\
&\le C \,  \Big\| \left(\frac{\partial u_\e}{\partial \nu_\e}\right)_+\Big\|_2 
+ C \, r_0^{-1/2} \| \nabla u_\e\|_{L^2(\Omega_-)}  + C \, r_0^{-1/2} \| \nabla u_\e\|_{L^2(\Omega_+)}, 
\endaligned
$$
where we also used the fact $(\nabla_{\tan} u_\e)_+ = (\nabla_{\tan}  u_\e)_-$ on $\partial\Omega$.
This, together with  Lemma \ref{lemma-3.20}, gives 
$$
\| f\|_2 \le C\, \Big \| \left(\frac{\partial u_\e}{\partial\nu_\e}\right)_+ \Big\|_2 
+ C\,  \| f\|_2 ^{1/2} \Big\| \left(\frac{\partial u_\e}{\partial\nu_\e}\right)_+ \Big\|_2 ^{1/2},
$$
which, by the Cauchy inequality with an $\e>0$, leads to the estimate (\ref{inverse-bound}).
}
\end{remark}

The remark above reduces Theorem \ref{inverse-theorem} to the following.

\begin{theorem}\label{Rellich-theorem}
Assume that  $A$ satisfies the conditions (\ref{elasticity}), (\ref{periodicity}) and (\ref{Holder}).
Let $\Omega$ be a bounded Lipschitz domain in $\rd$.
Let $u_\e$  be a weak solution of $\mathcal{L}_\e (u_\e)=0$ in $\Omega$
with the properties that $(\nabla u_\e)^*\in L^2(\partial\Omega)$ and
$\nabla u_\e$ exists n.t. on $\partial\Omega$.
Then
\begin{equation}\label{Rellich-estimate-000}
\int_{\partial\Omega} |\nabla u_\e|^2\, d\sigma
\le C \int_{\partial\Omega} \Big| \frac{\partial u_\e}{\partial\nu_\e} \Big|^2  d\sigma
+\frac{C}{r_0} \int_\Omega |\nabla u_\e|^2\, dx,
\end{equation}
where $r_0=\text{\rm diam}(\Omega)$ and
$C$ depends only on $\kappa_1$, $\kappa_2$, $(\sigma, M)$ and
the Lipschitz character of $\Omega$.
\end{theorem}

\begin{proof}[\bf Theorem \ref{Rellich-theorem} $\Longrightarrow$ Theorem \ref{inverse-theorem}]

It follows from Remark \ref{remark-3.20} that the Rellich estimate (\ref{Rellich-estimate-000})
implies  the inequality (\ref{inverse-bound}). 
To show $(1/2)I +\mathcal{K}_\e$ is invertible on $L^2_{\mathcal{R}} (\partial\Omega)$,
we apply a continuity argument by considering the matrix
$$
A^t (y) =t A(y)  + (1-t) A^0,
$$
where $t\in [0,1]$ and $A^0$ is a constant matrix with entries 
$\delta_{ij} \delta_{\alpha\beta} +\delta_{i\alpha}\delta_{j\beta} +\delta_{i\beta} \delta_{j\alpha}$.
Let $\mathcal{L}_\e^t =-\text{\rm div} (A^t(x/\e)\nabla )$ and $(1/2)I +\mathcal{K}_\e^t$ be 
the corresponding operator, associated with $\mathcal{L}_\e^t$.
Note that $A^t$ satisfies conditions (\ref{elasticity}), (\ref{periodicity}) and (\ref{Holder})
with constants depending only on $\kappa_1$, $\kappa_2$, $(\sigma, M)$.
It follows from Theorem \ref{Rellich-theorem} that the estimate (\ref{Rellich-estimate-000})
holds for solutions of $\mathcal{L}_\e^t (u_\e)=0$ in $\Omega$,
with constant $C$ depending only on $\kappa_1$, $\kappa_2$, $(\sigma, M)$ and
the Lipschitz character of $\Omega$.
Consequently, by Remark \ref{remark-3.20},
we obtain
\begin{equation}\label{3.100}
\| f\|_2 \le C \| \big( (1/2) I +\mathcal{K}^t_\e \big) f\|_2 \quad
 \text{ for any } f\in L^2_{\mathcal{R}}(\partial\Omega),
\end{equation}
where $C$ depends only on $\kappa_1$, $\kappa_2$, $(\sigma, M)$ and
the Lipschitz character of $\Omega$.
Observe that 
$$
\| \mathcal{K}_\e^ t - \mathcal{K}_\e^s \|_{L^2(\partial\Omega) \to L^2(\partial\Omega)}
\to 0 \quad \text{ as } t \to s
$$
\cite{KS2}.
Since $(1/2)I +\mathcal{K}_\e^t $ is invertible on $L^2_{\mathcal{R}}(\partial\Omega)$ 
for $t=0$ \cite{DKV-1988},
it follows from (\ref{3.100}) that $(1/2)I +\mathcal{K}_\e^t $ is invertible on $L^2_{\mathcal{R}}(\partial\Omega)$ 
for any $t\in [0,1]$.  In particular,
the case $t=1$ gives Theorem \ref{inverse-theorem}.
\end{proof}

Note that by Remark \ref{remark-3.0}, 
Theorem \ref{main-theorem-2} $\Longrightarrow$ Theorem \ref{Rellich-theorem}.
Thus we have proved that for each $\e>0$ and  each Lipschitz domain $\Omega$,
$$
\text{\rm Theorem \ref{main-theorem-2}}
\Longleftrightarrow
\text{\rm Theorem \ref{inverse-theorem}}
\Longleftrightarrow
\text{\rm Theorem \ref{Rellich-theorem}}
$$

The rest of this paper is devoted to the proof of Theorem \ref{Rellich-theorem}.



\section{Nontangential maximal function and square function estimates}
\setcounter{equation}{0}

In this section we establish some nontangential maximal function and square function estimates for
weak solutions of Dirichlet problem 
\begin{equation}\label{D-3.1}
\left\{
\aligned
-\text{\rm div} (A(x)\nabla u)  & =F_0 +\text{\rm div}(F_1) & \quad & \text{ in } \Omega,\\
u& =0 &\quad & \text{ on } \partial \Omega,
\endaligned
\right.
\end{equation}
where $\Omega$ is a bounded Lipschitz domain in $\rd$.
These auxiliary estimates will be used in the next section to prove the boundary Korn inequality.
Throughout the section we will assume that $A=(a_{ij}^{\alpha\beta} (x))$ 
with $1\le i, j \le d$ and $1\le \alpha, \beta \le m$ satisfies the ellipticity condition (\ref{ellipticity}) and
the H\"older continuity condition (\ref{Holder}).
We also need the symmetry condition $A^*=A$, i.e.,
 $a_{ij}^{\alpha\beta}=a_{ji}^{\beta\alpha}$.
The estimates we will prove are local.
Thus, without loss of generality, we further assume that
 $\Omega\subset B(0,1/4)$, diam$(\Omega)\approx 1$, and $A$ is 1-periodic.
 
 \begin{lemma}\label{lemma-3.1}
 Let $u$ be a weak solution of (\ref{D-3.1}) with $F_0=0$.
 Then for any $0\le \sigma_1<\sigma_2 \le 1$,
 \begin{equation}\label{3.1-0}
 \int_\Omega |\nabla u|^2 \big[\delta (x)\big]^{\sigma_2}\, dx \le C
 \int_\Omega |F_1|^2  \big[ \delta (x)\big]^{\sigma_1}\, dx,
 \end{equation}
 where $\delta (x)=\text{dist}(x, \partial\Omega)$ and
 $C$ depends only on $\mu$, $\sigma_1$, $\sigma_2$, $(\sigma, M)$ in (\ref{Holder})
  and the Lipschitz character of $\Omega$.
 \end{lemma}
 
 \begin{proof}
 This was proved in \cite[pp.373-374]{Shen-2015} under the assumption that
 $A$ satisfies (\ref{ellipticity}), (\ref{Holder}), and $A^*=A$.
 \end{proof}
 
 \begin{lemma}\label{lemma-3.2}
 Let $u$ be a weak solution of (\ref{D-3.1}) with $F_1=0$.
 Then for any $0\le \sigma_1<\sigma_2\le 1$,
 \begin{equation}\label{3.2-0}
 \int_\Omega |\nabla u|^2 \big[ \delta (x)\big]^{\sigma_2}\, dx \le C
 \int_\Omega |F_0|^2  \big[ \delta (x)\big]^{\sigma_1+2}\, dx,
 \end{equation}
 where 
 $C$ depends only on $\mu$, $\sigma_1$, $\sigma_2$, $(\sigma, M)$ in (\ref{Holder})
  and the Lipschitz character of $\Omega$.
  \end{lemma}
 
 \begin{proof}
 Let $0\le \sigma_1< \sigma_2\le 1$.
 Suppose that $-\text{\rm div} (A\nabla v)=\text{\rm div} (F)$ in $\Omega$ and
 $v=0$ on $\partial\Omega$.
 It follows from Lemma \ref{lemma-3.1} by duality that
 \begin{equation}\label{3.2-1}
 \int_\Omega |\nabla v (x) |^2 \frac{dx}{ \big[\delta (x) \big]^{\sigma_1}}
 \le C \int_\Omega |F(x)|^2 \frac{dx}{\big[ \delta (x)\big]^{\sigma_2}}.
 \end{equation}
 This, together with Hardy's inequality (see e.g. \cite{Stein}),
 \begin{equation}\label{Hardy}
 \int_\Omega |u(x)|^2 \frac{dx}{\big[ \delta (x)\big]^{\alpha}}
 \le C_\alpha  \int_\Omega |\nabla  u(x)|^2 \frac{dx}{\big[\delta (x)\big]^{\alpha -2}},
 \end{equation}
 which holds for $\alpha>1$ and $u\in H^1_0(\Omega)$, 
 gives
 \begin{equation}\label{3.2-3}
  \int_\Omega |v(x)|^2 \frac{dx}{\big[ \delta (x)\big]^{\sigma_1 +2}}
  \le C \int_\Omega |F(x)|^2 \frac{dx}{\big[ \delta (x)\big]^{\sigma_2}}.
 \end{equation}
 The estimate (\ref{3.2-0}) follows from (\ref{3.2-3}) by a duality argument.
 \end{proof}
 
By combining Lemmas \ref{lemma-3.1} and \ref{lemma-3.2} we obtain the following.

\begin{theorem}\label{theorem-3.1}
Assume that $A$ satisfies (\ref{ellipticity}), (\ref{Holder}), and $A^*=A$.
Let $\Omega$ be a bounded Lipschitz domain with
diam$(\Omega)\approx 1$.
Let $u$ be a weak solution of (\ref{D-3.1}). Then
\begin{equation}\label{square-function-estimate}
\int_\Omega |\nabla u|^2 \big[\delta (x)\big]^{\sigma_2}\, dx
\le C \int_\Omega | F_0|^2 \big[\delta (x)\big]^{\sigma_1 +2}\, dx
+C \int_\Omega |F_1|^2 \big[\delta (x)\big]^{\sigma_1}\, dx
\end{equation}
for any $0\le \sigma_1<\sigma_2 \le 1$, where 
$C$ depends only on $\mu$, $\sigma_1$, $\sigma_2$,
$(\sigma, M)$ in (\ref{Holder})
  and the Lipschitz character of $\Omega$.
\end{theorem}

If $\sigma_2=1$, estimate (\ref{square-function-estimate})
is a square function estimate.
Note that by dilation (\ref{square-function-estimate}) continues to hold for 
any Lipschitz domain in $\rd$, with constant $C$  depending also on diam$(\Omega)$.

The rest of this section is devoted to the estimate of the nontangential maximal function of $u$.

\begin{lemma}\label{lemma-3.4}
Let $u$ be a weak solution of (\ref{D-3.1}). Then
\begin{equation}\label{3.4-0}
|u(x)|\le C \average_{B(x, r)} |u|
+ C \int_{B(x, r)} \frac{|F_1 (y)|\, dy}{|x-y|^{d-1}}
+C \int_{B(x, r)} \frac{|F_0(y)|\, dy}{|x-y|^{d-2}}
\end{equation}
for any $x\in \Omega$, where $r=\delta(x)/4$.
\end{lemma}

\begin{proof}
This is a standard interior estimate and may be proved by using the fundamental solution
$\Gamma (x, y)$ for the operator $-\text{\rm div} (A\nabla)$ in $\rd$.
Indeed, since we may assume $A$ is 1-periodic, we have
$|\Gamma (x, y)|\le C |x-y|^{2-d}$ and 
$|\nabla_y \Gamma (x, y)|\le C |x-y|^{1-d}$ \cite{AL-1991}.
Fix $x_0\in \Omega$ and let
$$
v(x)=-\int_{B(x_0, r_0)} \nabla_y \Gamma (x, y) F_1 (y)\, dy
+\int_{B(x_0, r_0)} \Gamma(x, y) F_0 (y)\, dy, 
$$
where $r_0=\delta(x_0)/4$.
Then 
$
\text{\rm div}(A\nabla (u-v)) =0 \text{ in } B(x_0, r_0).
$
Hence, by the interior $L^\infty$ estimate,
$$
|u(x_0)- v(x_0)|\le C \average_{B(x_0, r_0)} | u-v|,
$$
from which the estimate (\ref{3.4-0}) follows easily.
\end{proof}

We introduce a modified nontangential maximal function.
For $z\in \partial\Omega$, define
\begin{equation}\label{M}
N(w) (z)=\sup_{x\in \gamma(z) } \average_{B(x, \delta(x)/4)} |w|,
\end{equation}
where $w$ is a function defined on $\Omega$ and
and $\gamma (z)=\{ x\in \Omega: |x-z|< C_0 \delta(x) \}$.

\begin{lemma}\label{lemma-3.5}
Let $N(w)$ be defined by (\ref{M}) and $0<\alpha<1$.
Then for any $w\in H_0^1(\Omega)$,
\begin{equation}\label{3.5-0}
\int_{\partial\Omega} | N(w)|^2\, d\sigma
\le C \int_\Omega |\nabla w(x)|^2 \big[\delta (x)\big]^{\alpha}\, dx,
\end{equation}
where $C$ depends only on $\alpha$ and the Lipschitz character of $\Omega$.
\end{lemma}

\begin{proof}
Let  $z\in \partial\Omega$ and $x\in \gamma (z)$.
Note that, if $\delta(x)$ is sufficiently small,
$$
\average_{B(x, \delta(x)/4)} |w|\,  dy
\le C \average_{B(z, C\delta (x))\cap \Omega} |w|\, dy
\le C\average_{B(z, C \delta(x))\cap \partial\Omega} |M_r (w)|\, d\sigma,
$$
where $M_r (w)$ denotes the radial maximal function of $w$ (see e.g. \cite{KLS2} for its definition).
It follows that
$$
N(w)(z) \le C M_{\partial\Omega} \big(M_r (w)\big) (z) + C \int_\Omega |w|,
$$
for any $z\in \partial\Omega$, where $M_{\partial\Omega}$ denotes the Hardy-Littlewood maximal function 
operator on $\partial\Omega$.
Hence, by the $L^2$ boundedness of $M_{\partial\Omega}$,
$$
\aligned
\int_{\partial\Omega} |N(w)|^2\, d\sigma
 & \le C \int_{\partial\Omega} | M_r(w)|^2 \, d\sigma + C \int_\Omega |w|^2\, dx\\
& \le C \int_\Omega |\nabla w(x)|^2 \big[\delta (x)\big]^{\alpha}\, dx
+ C \int_\Omega |w|^2\, dx,
\endaligned
$$
where the last inequality was proved in \cite{KLS2}.
This, together with Hardy's inequality (\ref{Hardy}), gives (\ref{3.5-0}).
\end{proof}

To handle the last two terms in the r.h.s. of (\ref{3.4-0}), we introduce another 
nontangential maximal function
\begin{equation}\label{M-1}
\mathcal{M} (w)(z)=\sup\Big\{ |w (x)|: \, x\in \widetilde{\gamma}(z) \Big\},
\end{equation}
where $z\in \partial\Omega$ and
$$
\widetilde{\gamma} (z) =\big\{ x\in \Omega: |x-z|< 10 C_0 \delta (x) \big\}\supset \gamma (z)
=\big\{ x\in \Omega: |x-z|< C_0 \delta (x)\big\}.
$$
Observe that if $x\in \gamma(z)$ and $y\in B(x, \delta(x)/4)$, then
$y\in \widetilde{\gamma}(z)$.

\begin{theorem}\label{theorem-3.2}
Assume that $A$ satisfies (\ref{ellipticity}), (\ref{Holder}), and $A^*=A$.
Let $\Omega$ be a bounded Lipschitz domain with
diam$(\Omega)\approx 1$.
Let $u$ be a weak solution of (\ref{D-3.1})
and $0<\alpha, \sigma_1<1$.
Then for any $0<t<1$,
\begin{equation}\label{3.6-0}
\aligned
\int_{\partial\Omega} |(u)^*|^2\, d\sigma
 &\le C t^{2\alpha} \int_{\partial\Omega} | \mathcal{M} (F_1 \delta^{1-\alpha}) |^2 \, d\sigma
+C t^{2\alpha} \int_{\partial\Omega} |\mathcal{M} (F_0\delta^{2-\alpha})|^2\, d\sigma\\
&\qquad\qquad
+ C_t \int_\Omega |F_0|^2 \big[\delta(x)\big]^{\sigma_1 +2}\, dx
+C_t \int_\Omega |F_1|^2 \big[\delta(x)\big]^{\sigma_1}\, dx,
\endaligned
\end{equation}
where $\delta =\delta (x)$ and
$C$ depends only on $\mu$, $\sigma_1$, $\alpha$, $(\sigma, M)$ and the Lipschitz
character of $\Omega$. The constant $C_t$ also depends on $t$.
\end{theorem}

\begin{proof}
Fix $t\in (0,1)$. Let $z\in \partial\Omega$ and $x\in \gamma (z)$.
Note that if $r=\frac{\delta(x)}{4}\le t$, then
$$
\int_{B(x, r)} \frac{| F_1 (y)|\, dy}{|x-y|^{d-1}}
\le C r^\alpha \mathcal{M} (F_1 \delta^{1-\alpha}) (z)
\le C t^\alpha \mathcal{M} (F_1 \delta^{1-\alpha}) (z).
$$
If $r\ge t$, we have 
$$
\aligned
\int_{B(x, r)} \frac{| F_1 (y)|\, dy}{|x-y|^{d-1}}
 &=\int_{B(x, t)} \frac{| F_1 (y)|\, dy}{|x-y|^{d-1}}
+\int_{B(x, r)\setminus B(x, t)}  \frac{| F_1 (y)|\, dy}{|x-y|^{d-1}}\\
&\le C t^\alpha \mathcal{M}(F_1 \delta^{1-\alpha}) (z)
+C_t \int_{\{ y\in \Omega: \, \delta (y)\ge ct \}} |F_1|.
\endaligned
$$
Similarly, we may show that 
$$
\int_{B(x, r)} \frac{|F_0 (y)|\, dy}{|x-y|^{d-2}}
\le C t^\alpha \mathcal{M}(F_0 \delta^{2-\alpha}) (z) +C_t \int_{\{ y\in \Omega: \, \delta (y)\ge ct \}} |F_0|.
$$
In view of Lemma \ref{lemma-3.4} we obtain 
$$
(u)^* \le C N (u) + C t^\alpha \Big\{ \mathcal{M} (F_0\delta^{2-\alpha})
+\mathcal{M}(F_1\delta^{1-\alpha}) \Big\}
+ C_t  \int_{\{ y\in \Omega: \, \delta (y)\ge ct \}} ( |F_0| +|F_1|).
$$
This, together with Lemma \ref{lemma-3.5} and Theorem \ref{theorem-3.1},
gives (\ref{3.6-0}).
\end{proof}


\section{Boundary Korn inequality}
\setcounter{equation}{0}

Throughout this section we assume that $A$ satisfies (\ref{elasticity}) and
(\ref{Holder}). Let $\Omega$ be a bounded Lipschitz domain such that
$\Omega\subset B(0,1/4)$.
We further assume that there exists $M_0>0$ such that
\begin{equation}\label{additional}
|\nabla A(x)|\le M_0  \big[ \delta (x)\big]^{\sigma-1} \quad
\text{ and } \quad |\nabla^2 A(x)|\le M_0  \big[ \delta (x)\big]^{\sigma-2},
\end{equation}
for any $x\in B(0,1/2)$.
The goal of this section is to prove the following.

\begin{theorem}\label{local-Korn}
Assume that $A$ and $\Omega$ satisfy the conditions stated above.
Let $u$ be a weak solution of $\mathcal{L}_1 (u)=0$ in $\Omega$ with
the properties that $(\nabla u)^*\in L^2(\partial\Omega)$ and $\nabla u$ exists n.t. on $\partial\Omega$.
Then
\begin{equation}\label{4.1}
\| \nabla u\|_2
\le C \, \| \nabla u + (\nabla u)^T\|_2
+C \, r_0^{-1/2} \| \nabla u\|_{L^2(\Omega)},
\end{equation}
where $r_0=\text{\rm diam}(\Omega)$ and
$C$ depends only on $\kappa_1$, $\kappa_2$, $(\sigma, M)$, $M_0$, 
and the Lipschitz character of $\Omega$.
\end{theorem}

Since $\Omega \subset B(0,1/4)$, by using a periodic extension of $A$, we may assume that $A$ is 1-periodic in $\rd$.
Also, by approximating $\Omega$ with a sequence of smooth domains from inside with uniform Lipschitz character \cite{Verchota-1984},
we may assume that $A\in C^2(\overline{\Omega})$ and $u\in C^3(\overline{\Omega})$.
The assumptions that $(\nabla u)^*\in L^2(\partial\Omega)$ and
$\nabla u$ exists n.t. on $\partial\Omega$ allow us to prove (\ref{4.1})
by a limiting argument.

We begin with a boundary Korn inequality for harmonic functions.

\begin{lemma}\label{harmonic-lemma}
Let $u=(u^1, u^2, \dots, u^d)$ be a solution of $\Delta u=0$ in $\Omega$.
Suppose that $(\nabla u)^*\in L^2(\partial\Omega)$.
Then the inequality (\ref{4.1}) holds with a constant $C$ depending only on
 the Lipschitz character of $\Omega$.
\end{lemma}

\begin{proof}
This was proved in \cite[p.804-805]{DKV-1988}.
Note that  for harmonic function $u$,
the condition $(\nabla u)^*\in L^2(\partial\Omega)$ implies  $\nabla u$ exists n.t. on $\partial\Omega$.
\end{proof}

To utilize Lemma \ref{harmonic-lemma}, we observe that by Proposition \ref{prop-2.2},
$\mathcal{L}_1 (u)=0$ in $\Omega$ implies
$\text{\rm div} (\widetilde{A} (x)\nabla u)=0$ in $\Omega$. It follows that
\begin{equation}\label{4.10}
\mu \Delta u^\alpha +\frac{\partial}{\partial x_i}
\left\{ b_{ij}^{\alpha\beta} (x) \frac{\partial u^\beta}{\partial x_j} \right\} =0 \quad \text{ in } \Omega,
\end{equation}
where $b_{ij}^{\alpha\beta}$ is given by (\ref{b}).
Let 
\begin{equation}\label{Gamma}
\Gamma (f) (x)=\int_\Omega \Gamma (x-y) f(y)\ dy,
\end{equation}
where $\Gamma (x)$ denotes the fundamental solution for $\Delta$ in $\rd$, with pole at the origin.
Then
$$
\Delta \big\{ \Gamma (f)\big\} =f \quad \text{ in } \Omega.
$$
This allows us to rewrite  (\ref{4.10}) as
\begin{equation}\label{4.11}
\Delta \left\{ \mu u^\alpha + \frac{\partial }{\partial x_i} \Gamma \left(
b_{ij}^{\alpha\beta} \frac{\partial u^\beta}{\partial x_j} \right) \right\} =0 \quad 
\text{ in } \Omega.
\end{equation}

\begin{lemma}\label{lemma-4.2}
Let $u\in C^3(\overline{\Omega}; \rd)$ be a solution of $\mathcal{L}_1 (u)=0$ in $\Omega$.
Then
\begin{equation}\label{4.2-1}
\| \nabla u\|_2
\le C \, \|\nabla u +(\nabla u)^T\|_2
+ C\,  \|\nabla^2 \Gamma (B\nabla u)\|_2
+Cr_0^{-1/2} \| \nabla u\|_{L^2(\Omega)},
\end{equation}
where $B=(b_{ij}^{\alpha\beta})$.
\end{lemma}

\begin{proof}
This follows directly from Lemma \ref{harmonic-lemma} and (\ref{4.11}).
We only need to observe that
$$
\|\nabla^2 \Gamma (B \nabla u)\|_{L^2(\Omega)} \le C\,  \| \nabla u\|_{L^2(\Omega)},
$$
which is a standard singular integral estimate.
\end{proof}

The estimate of the second term in the r.h.s. of (\ref{4.2-1}) is much involved.
We start with the following.

\begin{lemma}\label{lemma-4.3}
Let $\Omega$ be a bounded Lipschitz domain and $w\in C^1(\overline{\Omega})$.
Then
\begin{equation}\label{4.3-0}
\|\nabla^2 \Gamma (w)\|_2
\le C \left\{ \| (w)^*\|_2
+\left(\int_\Omega |\nabla w(x)|^2\delta (x)\, dx \right)^{1/2} \right\},
\end{equation}
where $C$ depends only on the Lipschitz character of $\Omega$.
\end{lemma}

\begin{proof}
Let $g\in L^2(\partial \Omega)$.
It follows by Fubini's Theorem that
$$
\aligned
\int_{\partial\Omega} \nabla_y^2 \Gamma (w) \cdot g(y)\, dy
&=\int_\Omega w(x) \cdot \left\{ \int_{\partial\Omega}
\nabla_y^2 \Gamma (y-x) g(y) dy\right\}\, dx\\
&=- \int_\Omega w(x) \cdot \left\{ \nabla_x \int_{\partial\Omega} \nabla_y \Gamma (y-x) g (y)\, dy \right\}\, dx.
\endaligned
$$ 
Since the function 
$$
v(x)=\int_{\partial\Omega} \nabla_y \Gamma (y-x) g (y)\, dy 
$$
is harmonic in $\Omega$, it follows from Dahlberg's bilinear estimate \cite{Dahlberg-1986}  that
$$
\aligned
\Big| \int_{\partial\Omega} \nabla_y^2 \Gamma (w) \cdot g(y)\, dy\Big|
 &\le C \, \| (v)^*\|_2
 \left\{ \| (w)^*\|_2
+\left(\int_\Omega |\nabla w(x)|^2\delta (x)\, dx \right)^{1/2} \right\}\\
& \le C\,  \| g\|_2
 \left\{ \| (w)^*\|_2
+\left(\int_\Omega |\nabla w(x)|^2\delta (x)\, dx \right)^{1/2} \right\},
\endaligned
$$
where we have used the well known estimate $\|(v)^*\|_2 \le C \| g\|_2$
\cite{Verchota-1984}. 
The desired estimate now follows by duality.
\end{proof}

Next, we factor  the matrix $B=(b_{ij}^{\alpha\beta}) $.
Let $m=d(d+1)/2$.

\begin{lemma}\label{lemma-4.4}
Let $B(x)=\big( b_{ij}^{\alpha\beta} (x)\big)$, where $b_{ij}^{\alpha\beta} $ is given by (\ref{b}).
Then there exists a matrix $Q (x)=\big(q_{i}^{t \alpha} (x) \big)$ with $1\le \alpha, i\le d$  and
$1\le t \le m$ such that $q_i^{t \alpha} =q_\alpha^{t i}$,
\begin{equation}\label{4.4-0}
b_{ij}^{\alpha\beta}= q_i^{t \alpha} q_j^{t \beta} \quad \text{ for any }
1\le \alpha, \beta, i, j \le d.
\end{equation}
Moreover, $|Q(x)|\le M_1$, $|Q(x)-Q(y)|\le M_1 |x-y|^\sigma$,
\begin{equation}\label{4.4-2}
|\nabla Q(x)|\le M_1 \big[ \delta (x)\big]^{\sigma -1} \quad \text{ and } \quad 
|\nabla^2 Q(x)|\le M_1 \big[ \delta (x)\big]^{\sigma -2}
\end{equation}
for any $x\in B(0,2)$, where $M_1$ depends only on $\kappa_1$, $\kappa_2$,
$(\sigma, M)$ and $M_0$ in (\ref{additional}).
\end{lemma}

\begin{proof}
We begin by fixing a constant matrix $(h_i^{t\alpha})$ with
$1\le \alpha, i \le d$ and $1\le t\le m$ such that $h_i^{t\alpha} =h_\alpha^{ti}$ and
$$
\big( h_i^{t\alpha} \xi_i^{\alpha}\big)_{1\le t\le m}= E (\xi)
\in \mathbb{R}^{m}
$$
for any $\xi =(\xi_i^\alpha)\in \rdd$,
where $E(\xi)$ is an enumeration of the lower triangular part of matrix $(1/2)(\xi +\xi^T)$, 
$$
E (\xi)= (e^t (\xi))_{1\le t\le m}= (\xi_1^1,  (\xi_2^1 +\xi_1^2)/2, \xi_2^2, \dots, \xi_d^d )
=\big( (\xi_j^\beta +\xi_\beta^j)/2 \big)_{1\le j\le d, 1\le \beta\le j}.
$$
Next we define the $m\times m$ symmetric matrix $G(x)=(g^{ts} (x))$ by the quadratic form
$$
\eta^t g^{ts} (x) \eta^s =\xi_i^\alpha b_{ij}^{\alpha\beta} (x)\xi_j^\beta,
$$
where $\eta  =(\eta^t)=E(\xi)\in \mathbb{R}^m$ and $\xi=(\xi^\alpha_i)\in \rdd$.
Using the fact that $E(\xi)=0$ if and only if $\xi^T=-\xi$ and that
$b_{ij}^{\alpha\beta}=b_{ji}^{\beta\alpha}=b_{\alpha j}^{i\beta}$,
it is easy to verify that $G(x)$ is well defined.
Moreover, in view of (\ref{b-ellipticity}), we obtain 
$$
\eta^t g^{ts} (x) \eta^s
\ge \frac{\kappa_1}{8} |\xi +\xi^T|^2 \ge \frac{\kappa_1}{8} |\eta|^2
$$
for any $\eta\in \mathbb{R}^m$. Thus $G(x)$ is an $m\times m$ 
symmetric, positive-definite matrix.
It follows that there exists a symmetric, positive-definite matrix $P=(p^{ts})$ such that
$G=P^2$. Moreover, since $G \ge {\kappa_1/8}$,
we see that $|P(x)|\le C$, $|P(x)-P(y)|\le C |x-y|^\sigma$,
$|\nabla P(x)|\le C \big[ \delta(x)\big]^{\sigma-1}$, 
and $|\nabla^2 P(x)|\le C \big[ \delta(x)\big]^{\sigma-2}$, 
where $C$ depends only on $\kappa_1$, $\kappa_2$,
$(\sigma, M)$, and $M_0$ in (\ref{additional}).
This man be proved by using the formula,
$$
P= G^{1/2}=\frac{1}{\pi} \int_0^\infty \lambda^{-1/2} G (\lambda I + G)^{-1}\, d\lambda
$$
and the resolvent identity, 
$$
(\lambda I + G(x))^{-1} -(\lambda I + G(y))^{-1}
=( \lambda I + G(x))^{-1} \big( G(y) -G(x) \big) ( \lambda I + G(y))^{-1}.
$$

Finally, since
$$
\xi_i^\alpha b_{ij}^{\alpha\beta} \xi_j^\beta
= h_i^{t\alpha} \xi_i^\alpha g^{ts} h_j^{s\beta} \xi_j^\beta
=\xi_i^\alpha h_i^{t \alpha}  p^{tk} p^{sk} h_j^{s\beta} \xi_j^\beta,
$$
for any $\xi=(\xi_i^\alpha) \in \rdd$,
we obtain 
$$
b_{ij}^{\alpha\beta} =h_i^{t\alpha} p^{tk}p^{sk} h_j^{s\beta} =q_i^{k\alpha} q_j^{k \beta},
$$
where we have set $q_i^{k\alpha} =q_i^{k \alpha}(x) =h_i^{t \alpha} p^{tk} (x)$.
The desired smoothness estimates for $Q=(q_i^{k\alpha}) $ follow from the same estimates for $P=(p^{ts})$.
\end{proof}

Let $v=(v^t)$, where $1\le t\le m$ and
\begin{equation}\label{4.21}
v^t =q_j^{t\beta}(x) \frac{\partial u^\beta}{\partial x_j}.
\end{equation}
Note that since $q_j^{t\beta} =q_\beta^{t j}$, we have $|v|\le C |\nabla u + (\nabla u)^T|$.

\begin{lemma}\label{lemma-4.5}
Let $u\in C^3(\overline{\Omega}; \rd)$ be a solution of $\mathcal{L}_1 (u)=0$ in $\Omega$.
Then
\begin{equation}\label{4.5-1}
\aligned
\| \nabla u\|_2
\le C  \|\nabla u +(\nabla u)^T\|_2
& + C \| (v)^*\|_2
+ C \left(\int_\Omega |\nabla v|^2 \delta (x)\, dx \right)^{1/2}\\
& +Cr_0^{-1/2} \| \nabla u\|_{L^2(\Omega)},
\endaligned
\end{equation}
where $v$ is given by (\ref{4.21}).
\end{lemma}

\begin{proof}
In view of Lemma \ref{lemma-4.2} it suffices to handle the term 
$\| \nabla^2 \Gamma (B \nabla u)\|_2$.
Let $w=B\nabla u=(w_i^\alpha)$, where 
$$
w_i^\alpha =b_{ij}^{\alpha\beta} \frac{\partial u^\beta}{\partial x_j}
=q_i^{t\alpha} q_j^{t\beta} \frac{\partial u^\beta}{\partial x_j}
=q_i^{t\alpha} v^t
$$
and we have used Lemma \ref{lemma-4.4}.
It follows by Lemma \ref{lemma-4.3} that
$$
\aligned
& \| \nabla^2 \Gamma (B \nabla u)\|_2\\
&\le C \| (w)^*\|_2
+ C \left(\int_\Omega |\nabla w (x)|^2 \delta (x)\, dx \right)^{1/2}\\
&\le C \| (v)^*\|_{L^2(\partial\Omega)}
+C \left(\int_\Omega |\nabla v|^2 \delta (x)\, dx \right)^{1/2}
+ C \left(\int_\Omega |v|^2 |\nabla Q|^2 \delta (x)\, dx \right)^{1/2}.
\endaligned
$$
Since $|\nabla Q(x)|\le C \big[\delta(x)\big]^{\sigma-1}$ for $x\in \Omega$, 
 we see that $|\nabla Q (x) |^2 \delta (x)\, dx$ is a Carleson measure on $\Omega$.
It follows that
$$
 \left(\int_\Omega | v|^2 |\nabla Q(x)|^2 \delta (x) dx \right)^{1/2}
\le C \| (v)^*\|_2
$$ 
This completes the proof.
\end{proof}

We are now in a position to give the proof of Theorem \ref{local-Korn}.

\begin{proof}[\bf Proof of Theorem \ref{local-Korn}]
It remains to control the second and third terms in the r.h.s. of (\ref{4.5-1}) and
show that
\begin{equation}\label{goal-4.1}
\| (v)^*\|_2
+ \left(\int_\Omega |\nabla v|^2 \delta (x)\, dx \right)^{1/2}
\le C \|\nabla u + (\nabla u)^T\|_2
+ C r_0^{-1/2} \| \nabla u\|_{L^2(\Omega)},
\end{equation}
where $v=Q (x)\nabla u(x)$ and $Q(x)$ is given by Lemma \ref{lemma-4.4}.
To this end we  rewrite the equation (\ref{4.10}) as
\begin{equation}\label{4.20}
\mu \Delta u^\alpha +\frac{\partial}{\partial x_i} \left\{ q_i^{t\alpha} q_j^{t\beta} \frac{\partial u^\beta}{\partial x_j}\right\}
=0 \quad \text{ in } \Omega.
\end{equation}
We now differentiate (\ref{4.20}) in $x_\ell$ and then multiply the resulting equation by
$q_\ell^{s\alpha}$ to obtain 
$$
\mu q_\ell^{s\alpha}  \Delta \frac{\partial u^\alpha}{\partial x_\ell}
+ q_\ell^{s\alpha} \frac{\partial^2}{\partial x_i \partial x_\ell} \Big\{ q_i^{t\alpha} v^t \Big\}=0
\quad \text{ in } \Omega.
$$
It follows that
$$
\mu \Delta v^s -\mu \Big[ \Delta, q_\ell^{s\alpha}\Big] \frac{\partial u^\alpha}{\partial x_\ell}
+q_\ell^{s\alpha} \frac{\partial}{\partial x_\ell} \left\{ \frac{\partial q_i^{t\alpha}}{\partial x_i} v^t \right\}
+ q_\ell^{s\alpha}\frac{\partial}{\partial x_\ell} \left\{ q_i^{t\alpha} \frac{\partial v^t}{\partial x_i} \right\}=0,
$$
where $[S,T]=ST-TS$ denotes the commutator of operators $S$ and $T$.
This gives
\begin{equation}\label{4.22}
\mu \Delta v^s +\frac{\partial }{\partial x_\ell} \left\{ q_\ell^{s\alpha }q_i^{t\alpha}\frac{\partial v^t}{\partial x_i} \right\}
=F^s \quad \text{ in } \Omega
\end{equation}
for $1\le s\le m$, where
\begin{equation}\label{4.23}
F^s= \mu \Big[ \Delta, q_\ell^{s\alpha}\Big] \frac{\partial u^\alpha}{\partial x_\ell}
-q_\ell^{s\alpha} \frac{\partial}{\partial x_\ell} \left\{ \frac{\partial q_i^{t\alpha}}{\partial x_i} v^t \right\}
+\frac{\partial q_\ell^{s\alpha}}{\partial x_\ell} \cdot q_i^{t\alpha} \frac{\partial v^t}{\partial x_i}.
\end{equation}
Using the product rule as well as the symmetry $q_\ell^{s\alpha}=q_\alpha^{s\ell}$,
a computation shows that
\begin{equation}\label{F-0}
F=\text{\rm div} (F_1) +F_0 \quad \text{ in  }\Omega,
\end{equation}
 where $F_0$ and $F_1$ satisfy the estimates
\begin{equation}\label{F}
\aligned
|F_0| & \le C \Big\{ |\nabla^2 Q| + |\nabla Q|^2 \Big\} |\nabla u +(\nabla u)^T|,\\
|F_1| & \le C |\nabla Q| | \nabla u +(\nabla u)^T|.
\endaligned
\end{equation}
We now decompose $v =  v_1+ v_2$ in $\Omega$, where
\begin{equation}\label{v-1}
\left\{ \aligned
\mu \Delta v_{1}^s +\frac{\partial }{\partial x_\ell} \left\{ q_\ell^{s\alpha }q_i^{t\alpha}\frac{\partial v_{1}^t}{\partial x_i} \right\}
 &= (\text{\rm div} (F_1))^s  +F_0^s &\quad & \text{ in } \Omega,\\
 v_{1}  & =   0& \quad & \text{ on } \partial\Omega,
 \endaligned
 \right.
 \end{equation}
 and
 \begin{equation}\label{v-2}
\left\{ \aligned
\mu \Delta v_2^s +\frac{\partial }{\partial x_\ell} \left\{ q_\ell^{s\alpha }q_i^{t\alpha}\frac{\partial v_{2}^t}{\partial x_i} \right\}
 &=0 &\quad & \text{ in } \Omega,\\
 v_2  & =   v& \quad & \text{ on } \partial\Omega.
 \endaligned
 \right.
 \end{equation}
To proceed, we  point out that the $m\times d$ system in (\ref{v-1}) and (\ref{v-2})
  is an elliptic system in divergence form
with coefficients 
$$
\mu \delta_{st}\delta_{\ell i} +q_\ell^{s\alpha} q_i^{t\alpha},
$$
which are very strongly elliptic, symmetric, and H\"older continuous.
As a result, it follows (see e.g.  \cite{KS2}) that
 the solution $v_2$ of Dirichlet problem (\ref{v-2}) satisfies the estimate
\begin{equation}\label{v-2-estimate}
\| (v_2)^*\|_2
+\left(\int_\Omega |\nabla v_2|^2 \delta (x)\, dx \right)^{1/2}
\le C \| v\|_2\\
\le C \| \nabla u + (\nabla u)^T\|_2.
\end{equation}

Finally, to estimate $v_1$, we use Theorems \ref{theorem-3.1} and \ref{theorem-3.2}.
This gives
\begin{equation}\label{4.25}
\aligned
& \int_{\partial\Omega}  | (v_1)^*|^2 \, d\sigma
+\int_\Omega |\nabla v_1|^2 \delta (x)\, dx \\
&\le C_t \int_\Omega |F_0|^2 \big[\delta (x)\big]^{\sigma_1 +2}\, dx
+ C_t \int_\Omega |F_1|^2 \big[ \delta (x)\big]^{\sigma_1}\, dx\\
&\qquad\qquad
+Ct^{2\alpha} \int_{\partial\Omega} |\mathcal{M}(F_1 \delta^{1-\alpha})|^2 \, d\sigma
+C t^{2\alpha} \int_{\partial\Omega} |\mathcal{M}(F_0 \delta^{2-\alpha})|^2 \, d\sigma,
\endaligned
\end{equation}
for any $t\in (0,1)$,
where $0< \alpha, \sigma_1<1$.
Recall that $|\nabla Q(x)|\le C \big[\delta(x)\big]^{\sigma-1}$ and
$|\nabla^2 Q(x)|\le C \big[ \delta (x)\big]^{\sigma-2}$.
In view of (\ref{F}) we obtain 
\begin{equation}\label{4.26}
\aligned
|F_0(x)| &\le C \big[\delta (x)\big]^{\sigma-2} |\nabla u + (\nabla u)^T|,\\
|F_1(x)| &\le C \big[\delta (x)\big]^{\sigma-1} |\nabla u + (\nabla u)^T|
\endaligned
\end{equation}
for any $x\in \Omega$.
By choosing $\alpha=\sigma$ and $\sigma_1 =1-\sigma$ we may deduce from (\ref{4.25}) and (\ref{4.26}) that
\begin{equation}\label{4.27}
\aligned
 &
  \int_{\partial\Omega}  | (v_1)^*|^2 \, d\sigma
+\int_\Omega |\nabla v_1|^2 \delta (x)\, dx \\
&\le C_t \int_\Omega |\nabla u +(\nabla u)^T|^2\, dx
+Ct^{2\sigma} \int_{\partial\Omega} |\mathcal{M} (\nabla u)|^2\, d\sigma\\
& \le C_t \int_\Omega |\nabla u +(\nabla u)^T|^2\, dx
+Ct^{2\sigma} \int_{\partial\Omega} |\nabla u |^2\, d\sigma
+Ct^{2\sigma} \int_{\partial\Omega} | u |^2\, d\sigma,
\endaligned
\end{equation}
where we have used (\ref{R-estimate-1}) for the last inequality.
This, together with (\ref{v-2-estimate}) and Lemma \ref{lemma-4.5},
gives
$$
\|\nabla u\|_2
  \le C \|\nabla u +(\nabla u)^T\|_2
+ C_t  \|\nabla u\|_{L^2(\Omega)}
+ C t^\sigma \|\nabla u\|_2
+  C t^{\sigma} \| u\|_2
$$
for any $t\in (0,1)$.
We now choose $t$ so small that $Ct^{\sigma}\le (1/2)$.
It follows that
$$
\aligned
\|\nabla u\|_2
   & \le C \|\nabla u +(\nabla u)^T\|_2
+ C  \|\nabla u\|_{L^2(\Omega)}
+  C \| u\|_2\\
& \le C \|\nabla u +(\nabla u)^T\|_2
+ C  \|\nabla u\|_{L^2(\Omega)}
+  C \| u\|_{L^2(\Omega)},
\endaligned
$$
where we have used the trace inequality (\ref{trace}) for the last step.
The proof is completed  by subtracting a constant from $u$ and using Poincar\'e inequality.
\end{proof}



\section{Rellich estimates for small scales}
\setcounter{equation}{0}

In this section we establish the Rellich estimate
\begin{equation}\label{Rellich}
\|\nabla u\|_2
\le C \Big \|\frac{\partial u}{\partial \nu}\Big\|_2
+ C r_0^{-1/2}  \| \nabla u\|_{L^2(\Omega)},
\end{equation}
for solutions of $\mathcal{L} (u)=-\text{\rm div}(A\nabla u)=0$ in a Lipschitz domain  $\Omega$,
and solve the $L^2$ Neumann problem in $\Omega$ by the method layer potentials.
This will be done under the assumptions that $A$ satisfies (\ref{elasticity}), (\ref{periodicity})
(\ref{Holder}) and  (\ref{additional}).
The extra assumption (\ref{additional}) will be eliminated in the next section by an approximation scheme.

Throughout this section $\Omega$ is a Lipschitz domain with
$r_0=\text{\rm diam}(\Omega) =1/4$.

\begin{lemma}\label{lemma-5.1}
Assume that $\mathcal{L}(u) =0$ in $\Omega$, $(\nabla u)^*\in L^2(\partial\Omega)$,
and $\nabla u $ exists n.t. on $\partial\Omega$.
Then
\begin{equation}\label{5.1-0}
\int_{\partial\Omega} |\nabla u + (\nabla u)^T |^2\, d\sigma
\le C \int_{\partial\Omega}  |\nabla u| \Big|\frac{\partial u}{\partial \nu} \Big|\, d\sigma
+C \int_\Omega (|\nabla A| +1) |\nabla u|^2\, dx,
\end{equation}
where $C$ depends only on $\kappa_1$, $\kappa_2$, and
the Lipschitz character of $\Omega$.
\end{lemma}

\begin{proof}
As in the case of constant coefficients \cite{DKV-1988},
the estimate (\ref{5.1-0}) follows from the so-called Rellich identity
\begin{equation}\label{5.1-1}
\int_{\partial\Omega} \langle \mathbf{h}, n\rangle a_{ij}^{\alpha\beta}
\frac{\partial u^\alpha}{\partial x_i} \cdot \frac{\partial u^\beta}{\partial x_j} \, d\sigma
=2\int_{\partial\Omega} \langle \mathbf{h},  \nabla u^\alpha \rangle \left(\frac{\partial u}{\partial \nu}\right)^\alpha
d\sigma +I,
\end{equation}
where $\mathbf{h}\in C_0^1(\rd, \rd)$,
$$
|I|\le C \int_\Omega \Big\{ |\nabla \mathbf{h}| +|\mathbf{h}| |\nabla A| \Big\} |\nabla u|^2\, dx,
$$
and $C$ depends only on $\kappa_1$ and $\kappa_2$.
The identity (\ref{5.1-1}) is proved by using integration by parts and the symmetry condition $a_{ij}^{\alpha\beta}
=a_{ji}^{\beta\alpha}$.
Since
\begin{equation}\label{5.1-2}
\frac{\kappa_1}{4}
|\nabla u +(\nabla u)^T|^2
\le a_{ij}^{\alpha\beta}
\frac{\partial u^\alpha}{\partial x_i} \cdot \frac{\partial u^\beta}{\partial x_j},
\end{equation}
the estimate (\ref{5.1-0}) follows by choosing $\mathbf{h}$ so that
$\langle \mathbf{h}, n \rangle \ge c>0$ on $\partial\Omega$.
\end{proof}

\begin{lemma}\label{lemma-5.2}
Assume that $\mathcal{L}(u) =0$ in $\Omega$, $(\nabla u)^*\in L^2(\partial\Omega)$,
and $\nabla u$ exists n.t. on $\partial\Omega$.
Then for any $t\in (0, 1)$,
\begin{equation}\label{5.2-0}
\|\nabla u\|_2 
\le C \Big\|\frac{\partial u}{\partial \nu} \Big\|_2
+ Ct^{\sigma/2}  \| (\nabla u)^*\|_2
+ C_t \|\nabla u\|_{L^2(\Omega)},
\end{equation}
where $C$ depends only on $\kappa_1$, $\kappa_2$, $(M, \sigma)$, and the Lipschitz
character of $\Omega$.
The constant $C_t$ also depends on $t$.
\end{lemma}

\begin{proof}
It follows by  the boundary Korn inequality in Theorem \ref{local-Korn} and (\ref{5.1-0}) that
\begin{equation}\label{5.2-1}
\| \nabla u \|_2 
\le C \Big\|\frac{\partial u}{\partial \nu} \Big\|_2
+ C \|\nabla u\|_{L^2(\Omega)}
+ C\| |\nabla A|^{1/2} \nabla u\|_{L^2(\Omega)},
\end{equation}
where we have also used the Cauchy inequality with an $\e$.
To estimate the last term in (\ref{5.2-1}), we use the assumption that
$|\nabla A (x)|\le C \big[\delta (x)\big]^{\sigma -1}$ for $x\in \Omega$.
This gives
\begin{equation}\label{5.2-2}
\aligned
\int_\Omega |\nabla A||\nabla u|^2\, dx
 &=\int_{\Omega_t} |\nabla A| |\nabla u|^2\, dx +\int_{\Omega\setminus \Omega_t} |\nabla A|
|\nabla u|^2\, dx\\
& \le C \int_{\Omega_t} \big[\delta (x)\big]^{\sigma-1} |\nabla u|^2\, dx
+C_t \int_\Omega |\nabla u|^2\, dx\\
&\le C t^\sigma \int_{\partial\Omega} |(\nabla u)^*|^2\, d\sigma
+ C_t \int_\Omega |\nabla u|^2\, dx,
\endaligned
\end{equation}
where $\Omega_t =\big\{ x\in \Omega: \ \delta (x) < t\big\}$.
The estimate (\ref{5.2-0}) now follows from (\ref{5.2-1}) and (\ref{5.2-2}).
\end{proof}

Let $(1/2) I +\mathcal{K}$ be the operator in (\ref{co-normal-limit}) with $\e=1$. 

\begin{theorem}\label{theorem-5.1}
Suppose that $A$ satisfies conditions (\ref{elasticity}), (\ref{periodicity}),
and (\ref{Holder}).
Let $\Omega$ be a bounded Lipschitz domain in $\rd$ such that
$\Omega \subset B(0,1)$ and diam$(\Omega)=(1/4)$.
 Also assume that $A$ satisfies the condition (\ref{additional}).
 Then $(1/2)I +\mathcal{K}: L^2_{\mathcal{R}} (\partial\Omega) \to L^2_{\mathcal{R}} (\partial\Omega)$
 is invertible and
 \begin{equation}\label{5.3-0}
 \| f\|_2 \le C \, \| \big( (1/2)I +\mathcal{K}\big) f\|_2,
 \end{equation}
 where $C$ depends only on $\kappa_1$, $\kappa_2$, $(M, \sigma)$, $M_0$ in (\ref{additional}),
 and the Lipschitz character of $\Omega$.
\end{theorem}

\begin{proof}
Let $u=\mathcal{S}(f)$ be the single layer potential for the operator 
$\mathcal{L}=\mathcal{L}_1$, where $f\in L^2_{\mathcal{R}} (\partial\Omega)$.
By  the jump relation (\ref{jump}) it follows that
$$
\|f\|_2 \le C\,  \| \left(\frac{\partial u}{\partial\nu}\right)_+\|_2 + C \| (\nabla u)_-\|_2.
$$
To deal with  $\|(\nabla u)_-\|_2$, we use Theorem \ref{Rellich-theorem-1} to obtain 
$$
\aligned
\| (\nabla u)_- \|_2  & \le  C\,  \| (\nabla_{\tan} u)_-\|_2
+ C\,  \| \nabla u\|_{L^2(\Omega_-)}\\
&\le C\,  \|(\nabla u)_+\|_2 + C\, \| \nabla u\|_{L^2(\Omega_-)},
\endaligned
$$
where we have used the fact (\ref{nt-limit}).
This, together with (\ref{5.2-0}), leads to
$$
\| f\|_2 \le C  \, \| \left(\frac{\partial u}{\partial\nu}\right)_+\|_2 
+ C t^{\sigma/2} \| (\nabla u)^*\|_2 + C_t \|\nabla u\|_{L^2(\Omega)}
+ C \|\nabla u\|_{L^2(\Omega_-)}
$$
for any $t\in (0,1)$. Since $\|(\nabla u)^*\|_2 \le C \| f\|_2$,
by choosing $t$ sufficiently small, we obtain 
\begin{equation}\label{5.3-10}
\aligned
\| f\|_2  & \le C  \| \left(\frac{\partial u}{\partial\nu}\right)_+\|_2 
+C \|\nabla u\|_{L^2(\Omega)}
+ C \|\nabla u\|_{L^2(\Omega_-)}\\
& \le C  \| \left(\frac{\partial u}{\partial\nu}\right)_+\|_2 
+ C \| f\|^{1/2} _2 \| \left(\frac{\partial u}{\partial\nu}\right)_+\|^{1/2}_2, 
\endaligned
\end{equation}
where we have used Lemma  \ref{lemma-3.20} for the last inequality.
By using the Cauchy inequality with an $\e>0$, this gives
\begin{equation}\label{5.3-20}
\| f\|_2 \le C\,   \| \left(\frac{\partial u}{\partial\nu}\right)_+\|_2 
= C\,   \| \big( (1/2)I +\mathcal{K}\big) f\|_2
\end{equation}
for any $f\in L^2_{\mathcal{R}} (\partial\Omega)$.

Finally, we recall that in the case of an isotropic body, where
$
a_{ij}^{\alpha\beta} =\lambda \delta_{i\alpha}\delta_{j\beta}
+\mu \delta_{ij}\delta_{\alpha\beta} 
+\mu \delta_{i \beta}\delta_{j \alpha},
$
the invertibility  of $(1/2)I + \mathcal{K}$ on $L^2_{\mathcal{R}}(\partial\Omega)$
was established in \cite{DKV-1988}.
By a simple continuity argument, this, together with the estimate (\ref{5.3-20}),
yields the invertibility  of $(1/2)I + \mathcal{K}$ on $L^2_{\mathcal{R}}(\partial\Omega)$
under the conditions stated in Theorem \ref{theorem-5.1}.
This completes the proof.
\end{proof}


\section{$L^2$ Neumann problems for small scales}
\setcounter{equation}{0}

In this section we use an approximate scheme \cite{KS2} to get rid of the extra assumption 
(\ref{additional}) in Theorem \ref{theorem-5.1}.
This allows us to solve the $L^2$ Neumann problem (\ref{NP})
in the small-scale case where $\e=1$ and diam$(\Omega)\le 1/4$.

\begin{theorem}\label{theorem-7.1}
Suppose that $A$ satisfies conditions (\ref{elasticity}), (\ref{periodicity}),
and (\ref{Holder}).
Let $\Omega$ be a bounded Lipschitz domain in $\rd$ such that
 diam$(\Omega)\le (1/4) $.
 Then $(1/2)I +\mathcal{K}: L^2_{\mathcal{R}} (\partial\Omega) \to L^2_{\mathcal{R}} (\partial\Omega)$
 is invertible and
 \begin{equation}\label{7.1-0}
 \| f\|_2 \le C\,  \| \big( (1/2)I +\mathcal{K}\big) f\|_2\quad \text{ for any } f\in L^2_\mathcal{R}(\partial\Omega),
 \end{equation}
 where $C$ depends only on $\kappa_1$, $\kappa_2$, $(\sigma, M)$, 
 and the Lipschitz character of $\Omega$.
\end{theorem}

\begin{remark}\label{remark-7.1}
{\rm
Theorem \ref{theorem-7.1} continues to hold without the periodicity condition (\ref{periodicity})
and for Lipschitz domains $\Omega$ with $r_0 =\text{\rm diam}(\Omega)\ge (1/4)$.
However, in this case, the constant $C$ will also depend on $r_0$.
This can be seen by a simple rescaling argument.
With the periodicity condition (\ref{periodicity}), the estimate (\ref{7.1-0})
holds with a constant $C$ independent of $r_0$.
This is equivalent to the same estimate for $\mathcal{K}_\e$ with constant $C$ independent of $\e$,
one of the main results of this paper.
}
\end{remark}

The reduction of Theorem \ref{theorem-7.1} from Theorem \ref{theorem-5.1} is similar to that in the case of
very strong ellipticity condition \cite{KS2}. We only provide an outline here.

Step 1. By translation and dilation  we may assume that $\Omega\subset B(0,1/4)$ and
diam $(\Omega)=(1/4)$.
We construct a coefficient matrix $\overline{A} (x) = (\overline{a}_{ij}^{\alpha\beta} (x))$,
with $1\le \alpha, \beta, i, j \le d$,
 in $\rd$ with the properties that 

(1) $\overline{A}=A$ on $\partial\Omega$;
 
 (2)  $\overline{A}$ satisfies the conditions (\ref{elasticity}), (\ref{periodicity}) and (\ref{Holder}),
 with possibly different constants depending only on $\kappa_1$, $\kappa_2$, $(\sigma, M)$ and the 
 Lipschitz character of $\Omega$;
 
 (3)  
 $\overline{A}$ satisfies the smoothness condition (\ref{additional}) with constant $M_0$ 
 depending on $(\sigma, M)$ and the Lipschitz character of $\Omega$.
 
\noindent  This is done as follows.
In $\Omega$ we let $\overline{a}_{ij}^{\alpha\beta}$ be the Poisson  extension 
of $A$, i.e.,
$$
\Delta \overline{a}_{ij}^{\alpha\beta}  =0 \quad \text{ in } \Omega \quad \text{ and } \quad 
\overline{a}_{ij}^{\alpha\beta}  =a_{ij}^{\alpha\beta}\quad \text{ on } \partial\Omega.
$$
On $[-1/2, 1/2]^d\setminus \overline{\Omega}$, we define
$\overline{a}_{ij}^{\alpha\beta}$ to be the harmonic function in $(-1/2, 1/2)^d\setminus \overline{\Omega}$
with boundary data $a_{ij}^{\alpha\beta}$ on $\partial\Omega$ and
$\delta_{i\alpha}\delta_{j\beta} +\delta_{ij}\delta_{\alpha\beta} +\delta_{i\beta}\delta_{j\alpha}$
on $\partial [-1/2, 1/2]^d$.
We then extend $\overline{A}$ to $\rd$ by periodicity.
The fact that $\overline{A}$ satisfies (\ref{elasticity}), (\ref{Holder}) and (\ref{additional})
follows from the maximum principle and well known estimates 
for harmonic functions in Lipschitz domains with H\"older continuous data.

Step 2. 
Let $\theta\in C_0^\infty (-1/2, 1/2)$ such that $0\le \theta\le 1$
and $\theta=1$ on $(-1/4, 1/4)$.
Define
\begin{equation}\label{7.9}
A^t (x) = \theta \left(\frac{\delta (x)}{t}\right) A(x) 
+\Big[ 1-\theta \left(\frac{\delta(x)}{t}\right) \Big] \overline{A} (x)
\end{equation}
for $x\in [-1/2, 1/2]^d$, where $t\in (0, 1/8)$ and $\overline{A}(x)$ is the matrix constructed in 
Step 1. Extend $A^t$ to $\rd$ by periodicity.
Let $\mathcal{K}_{A^t}$ denote the operator on $\partial\Omega$, associated with the conormal derivative 
of the single layer potential for $-\text{\rm div} \big(A^t (x)\nabla)$, as in Theorem \ref{layer-potential-theorem}.
Then
\begin{equation}\label{diff}
\| \mathcal{K}_{A^t} -\mathcal{K}_{\overline{A}} \|_{L^2(\partial\Omega)
\to L^2(\partial\Omega)}
\le C t^{\lambda_0},
\end{equation}
where $\lambda_0>0$ depends only on $\sigma$ and the Lipschitz character of $\Omega$.
We point out that this estimate was proved in \cite{KS2}
for operators satisfying the very strong  ellipticity condition (\ref{ellipticity}).
Since its proof only involves fundamental solutions, in view of
the equivalence of $\mathcal{L}_1$ and $\widetilde{\mathcal{L}}_1$,
it continues to hold for systems of elasticity.

Write
$$
(1/2) I +\mathcal{K}_{A^t}
=(1/2) I +\mathcal{K}_{\overline{A}} + (\mathcal{K}_{A^t}-\mathcal{K}_{\overline{A}}).
$$
Note that by Theorem \ref{theorem-5.1}, 
$$
\| \big((1/2) I + \mathcal{K}_{\overline{A}} \big)^{-1} \|_{L^2_{\mathcal{R}}(\partial\Omega)
\to L^2_{\mathcal{R}}(\partial\Omega)} \le C.
$$
It follows by (\ref{diff}) that there exist $t\in (0, 1/8)$ and $C>0$, depending only on $\kappa_1$, $\kappa_2$,
$(\sigma, M)$ and the Lipschitz character of $\Omega$, such that
\begin{equation}\label{7.10}
\| \big((1/2) I + \mathcal{K}_{A^t }\big)^{-1} \|_{L^2_{\mathcal{R}}(\partial\Omega)
\to L^2_{\mathcal{R}}(\partial\Omega)} \le C.
\end{equation}
Note that by (\ref{7.9}),
\begin{equation}\label{7.11}
A^t(x)= A(x) \quad \text{ if } \delta (x)\le (t/4).
\end{equation}

Step 3. It follows from (\ref{7.10}) that the $L^2$ Neumann problem (\ref{NP}) in $\Omega$
for the operator $\mathcal{L}^t =-\text{\rm div}(A^t(x)\nabla )$ is solvable and 
the estimate (\ref{N-estimate}) holds.
By Remark \ref{remark-3.0} this implies that if $\mathcal{L}^t (w)=0$ in $\Omega$,
$(\nabla w)^*\in L^2(\partial\Omega)$ and $\nabla w$ exists n.t. on $\partial\Omega$, then
\begin{equation}\label{7.12}
\int_{\partial\Omega}
|\nabla w|^2 \, d\sigma \le C \int_{\partial\Omega} \Big|\frac{\partial w}{\partial \nu}\Big|^2 d\sigma
+ C \, \Big| \int_\Omega \nabla w\, dx\Big|^2,
\end{equation}
where  $C$
depends only on $\kappa_1$, $\kappa_2$,
$(\sigma, M)$ and the Lipschitz character of $\Omega$.
Here we have used the fact that $A^t=A$ on $\partial\Omega$ and thus
the conormal derivative of $w$ on $\partial\Omega$ associated with $\mathcal{L}^t$
is the same as that associated with $\mathcal{L}=-\text{div}(A\nabla)$.

Now let $u$ be a weak solution of $\mathcal{L}(u)=0$ in $\Omega$ with the properties that
$(\nabla u)^*\in L^2(\partial\Omega)$ and $\nabla u$ exists n.t. on $\partial\Omega$.
Let $\varphi \in C_0^\infty(\rd)$ such that $0\le \varphi\le 1$,
$\varphi (x)=1$ on $\{ x\in \rd: \delta(x) \le (t/8) \}$, and
$\varphi (x)=0 $ on $\{ x\in \rd: \delta(x) \ge (t/6) \}$.
Let $\overline{u}=\varphi (u-E)$, where $E$ is the $L^1$ average of $u$ over $\Omega$.
Note that by (\ref{7.11}), we have
$\mathcal{L}^ t (u)=\mathcal{L} (u)=0$ on $\{ x\in \Omega: \delta (x) < (t/4) \}$.
It follows that
$$
\mathcal{L}^t (\overline{u})
=-\text{\rm div} (A^t \nabla \varphi \cdot (u-E)) -\nabla \varphi \cdot A^t \nabla u\quad \text{ in } \Omega.
$$
Let
$$
v(x)=\int_\Omega \nabla_y  \Gamma_t (x, y) 
A^t (y) \nabla \varphi  \cdot (u (y)-E)\, dy
-\int_\Omega \Gamma_t (x, y) \nabla \varphi  \cdot A^t (y)\nabla u (y)\, dy,
$$
where $\Gamma_ t(x, y)$ denotes the matrix of fundamental solutions for $\mathcal{L}^t$
in $\rd$. Then
$\mathcal{L}^t (\overline{u} -v)=0$ in $\Omega$ and
for any $x\in \Omega$ with $\delta (x)\le (t/16)$,
\begin{equation}\label{7.13}
\aligned
|\nabla v(x)| +| v(x)|
&\le C \int_\Omega | u-E|  +\int_\Omega |\nabla u|\\
&\le C \left(\int_\Omega |\nabla u|^2  \right)^{1/2},
\endaligned
\end{equation}
where we have used the observation $\nabla \varphi=0$ on 
$\{ x\in \Omega: \delta(x)\le (t/8)\}$.
This allows us to use the Rellich estimate (\ref{7.12}) for $w=\overline{u}-v$ and 
obtain 
$$
\aligned
\int_{\partial\Omega} |\nabla u|^2\, d\sigma
&\le 2 \int_{\partial\Omega} |\nabla w|^2\, d\sigma + 2\int_{\partial\Omega} |\nabla v|^2\, d\sigma\\
&\le C \int_{\partial\Omega} \Big|\frac{\partial u}{\partial \nu} \Big|^2d\sigma
+C \int_{\partial\Omega} |\nabla v|^2\, d\sigma
+C \, \Big|\int_\Omega \nabla u\Big|^2
+ C\,  \Big| \int_\Omega \nabla v\Big|^2\\
&\le C \int_{\partial\Omega} \Big|\frac{\partial u}{\partial \nu} \Big|^2d\sigma
+C\,   \Big|\int_\Omega \nabla u\Big|^2
+C \int_{\partial\Omega} |\nabla v|^2\, d\sigma
+ C \int_{\partial\Omega} | v|^2\, d\sigma,
\endaligned
$$
which, by (\ref{7.13}), yields the Rellich estimate for $u$,
\begin{equation}\label{7.14}
\int_{\partial\Omega} |\nabla u|^2\, d\sigma
\le C \int_{\partial\Omega} \Big|\frac{\partial u}{\partial \nu} \Big|^2d\sigma
+ C \int_\Omega |\nabla u|^2\, dx.
\end{equation}
Theorem \ref{theorem-7.1} now follows from
the equivalence of Theorem \ref{inverse-theorem} and Theorem \ref{Rellich-theorem}
for $\e=1$, proved in Section 3.

\begin{remark}\label{remark-7.2}
{\rm 
Let 
\begin{equation}\label{r-7.2-1}
\aligned
D (r, \psi) & = \Big\{ (x^\prime, x_d)\in \rd:
|x^\prime|<r \text{ and } \psi (x^\prime) < x_d < \psi (x^\prime) +10 (C_0+1)r \Big\},\\
\Delta (r, \psi) & = \Big\{ (x^\prime, \psi (x^\prime))\in \rd: |x^\prime|<r \Big\},
\endaligned
\end{equation}
where $\psi: \mathbb{R}^{d-1} \to \mathbb{R}$ is a Lipschitz function such that
$\psi (0)=0$ and $\|\nabla \psi\|_\infty \le C_0$.
Suppose that $\mathcal{L}_\e (u_\e)=0$ in $D(2r, \psi)$, where $0<r< c_0 \e$.
Let $v(x)= \e^{-1} u_\e (\e x)$. Then
$$
\mathcal{L}_1 (v) =0  \quad \text{ in } D({2r}/{\e}, \psi_\e),
$$
where $\psi_\e (x^\prime) =\e^{-1} \psi (\e x^\prime)$.
Under the conditions (\ref{elasticity}) and (\ref{Holder}), it follows from Theorem \ref{theorem-7.1}
that
\begin{equation}\label{r-7.2-2}
\int_{\Delta (tr/\e, \psi_\e)} |\nabla v|^2\, d\sigma
 \le C \int_{\partial D(tr/\e, \psi_\e)} \Big|\frac{\partial v}{\partial \nu} \Big|^2\, d\sigma
 +\frac{C\e}{r} \int_{D(tr/\e, \psi_\e)} |\nabla v|^2\, dx
 \end{equation}
for $1\le t\le 2$.
We point out that since $\psi_\e (0)=0$ and $\|\nabla \psi_\e\|_\infty \le C_0$,
the Lipschitz character of $D(tr/\e, \psi_\e)$ depends only on $C_0$.
As a result, the constant $C$ in (\ref{r-7.2-2}) depends only on
$\kappa_1$, $\kappa_2$, $(\sigma, M)$, $c_0$ and $C_0$.
By a change of variables we may deduce from (\ref{r-7.2-2}) that
\begin{equation}\label{r-7.2-3}
\aligned
\int_{\Delta (r, \psi)} |\nabla u_\e|^2\, d\sigma
\le C \int_{\Delta (2r, \psi)} \Big|\frac{\partial u_\e}{\partial \nu_\e} \Big|^2\, d\sigma
& +C \int_{\partial D(tr, \psi)\setminus \Delta(2r, \psi)} |\nabla u_\e|^2\, d\sigma\\
&+\frac{C}{r} \int_{D(2r, \psi)} |\nabla u_\e|^2\, dx.
\endaligned
\end{equation}
We now integrate both sides of (\ref{r-7.2-3}) with respect to $t$ over the interval $(1,2)$.
This gives
\begin{equation}\label{r-7.2-4}
\int_{\Delta (r, \psi)} |\nabla u_\e|^2\, d\sigma
\le C \int_{\Delta (2r, \psi)} \Big|\frac{\partial u_\e}{\partial \nu_\e} \Big|^2\, d\sigma
+\frac{C}{r} \int_{D(2r, \psi)} |\nabla u_\e|^2\, dx,
\end{equation}
where $0<r<c_0 \e$, $A$ satisfies (\ref{elasticity}) and (\ref{Holder}), and
$C$ depends only on
$\kappa_1$, $\kappa_2$, $(\sigma, M)$, $c_0$ and $C_0$.
}
\end{remark}


\section{Proof of Theorems \ref{main-theorem-2}, \ref{inverse-theorem} and \ref{Rellich-theorem}}
\setcounter{equation}{0}

As these three theorems are equivalent to each other (see Section 3),
it suffices to prove Theorem \ref{Rellich-theorem}.

We begin with a large-scale Rellich estimate.
Note that the smoothness condition (\ref{Holder}) is not needed.

\begin{theorem}\label{large-scale-Rellich-theorem}
Assume that $A$ satisfies conditions (\ref{elasticity}) and (\ref{periodicity}).
Let $\Omega$ be a bounded Lipschitz domain in $\rd$ and $g\in L^2_{\mathcal{R}}(\partial\Omega)$.
Let $u_\e\in H^1(\Omega; \rd)$ be a weak solution to the Neumann problem,
\begin{equation}\label{NP-8}
\mathcal{L}_\e (u_\e)=0  \quad \text{ in } \Omega
\quad \text{ and } \quad \frac{\partial u_\e}{\partial \nu_\e}=g 
\quad \text{ on } \partial\Omega.
\end{equation}
Suppose that $u_\e \perp \mathcal{R}$ in $L^2(\Omega; \rd)$.
Then for any $\e \le r< \text{\rm diam}(\Omega)$,
\begin{equation}\label{large-Rellich}
\frac{1}{r} \int_{\Omega_r} |\nabla u_\e |^2 \, dx \le C \int_{\partial\Omega} |g|^2\, d\sigma,
\end{equation}
where  $\Omega_r =\{ x\in \Omega: 
\delta(x)<r \}$ and $C$ depends only on 
$\kappa_1$, $\kappa_2$, and the Lipschitz character of $\Omega$.
\end{theorem}

\begin{proof}
This was proved in \cite[Theorem 1.2]{Shen-2016}.
\end{proof}

\begin{proof}[\bf Proof of Theorem \ref{Rellich-theorem}]
Suppose that $A$ satisfies (\ref{elasticity}), (\ref{periodicity}) and (\ref{Holder}).
Let $\Omega$ be a bounded Lipschitz domain in $\rd$ and $u_\e \in H^1(\Omega; \rd)$ a weak solution 
of $\mathcal{L}_\e (u_\e)=0$ in $\Omega$.
Assume that $(\nabla u_\e)^*\in L^2(\partial\Omega)$ and $\nabla u_\e$ exists n.t. on $\partial\Omega$.
We need to show that
\begin{equation}\label{Rellich-8}
\int_{\partial \Omega} |\nabla u_\e|^2\, d\sigma \le C \int_{\partial\Omega} |g|^2\, d\sigma
+\frac{C}{r_0} \int_\Omega |\nabla u_\e|^2\, dx,
\end{equation}
where $r_0=\text{\rm diam}(\Omega)$ and $g=\partial u_\e/{\partial \nu_\e}$.
We may assume that $0<\e< cr_0$.
The case $\e\ge cr_0$ follows easily from the local Rellich estimate in 
Section 7 for $\e=1$ by a simple rescaling argument.

First, we note that it suffices to prove (\ref{Rellich-8}) under the additional assumption 
that $ u_\e \perp \mathcal{R}$ in $L^2(\Omega; \rd)$.
Indeed, to handle the general case, we choose $\phi \in \mathcal{R}$ so that
$u_\e -\phi \perp \mathcal{R}$ in $L^2(\Omega; \rd)$. 
Then
\begin{equation}\label{8.1-1}
\aligned
\int_{\partial\Omega} |\nabla u_\e|^2\,  d \sigma
& \le 2 \int_{\partial\Omega} |\nabla u_\e -\nabla \phi |^2\,  d \sigma + 2\int_{\partial\Omega} |\nabla \phi |^2\,  d \sigma\\
 & \le C \int_{\partial\Omega} |g|^2\, d\sigma
+ \frac{C}{r_0} \int_\Omega |\nabla u_\e|^2\, dx
+\frac{C}{r_0} \int_\Omega |\nabla \phi |^2\, dx, 
\endaligned
\end{equation}
where we have used the fact that $\nabla \phi$ is constant.
Using the observation that $\phi$ is a linear function,
$$
\int_\Omega |\phi|^2\, dx  \le \int_\Omega |u_\e|^2\, dx  \quad \text{ and } \quad \int_\Omega u_\e\, dx
=\int_\Omega \phi \, dx ,
$$
as well as Poincar\'e inequality,
we obtain
$$
\aligned
\frac{1}{r_0} \int_\Omega |\nabla \phi|^2\, dx
 &\le \frac{C}{ r_0^3} \int_\Omega \Big| \phi -\average_\Omega \phi\Big|^2\, dx
\le \frac{C}{ r_0^3} \int_\Omega \Big| u_\e -\average_\Omega u_\e \Big|^2\, dx\\
&\le \frac{C}{r_0} \int_\Omega |\nabla u_\e|^2 \, dx,
\endaligned
$$
which, together with (\ref{8.1-1}), gives (\ref{Rellich-8}).

Next, fix $z\in \partial\Omega$.
It follows from the local Rellich estimate (\ref{r-7.2-4}) by a change of coordinate systems  that
\begin{equation}\label{8.1-2}
\int_{B(z, c \e)\cap \partial\Omega} |\nabla u_\e|^2\, d\sigma
\le C \int_{B(z, \e)\cap\partial\Omega} |g|^2\, d\sigma
+\frac{C}{\e} \int_{B(z,  \e)\cap \Omega} |\nabla u_\e|^2\, dx,
\end{equation}
where $C$ depends only on $\kappa_1$, $\kappa_2$, $(\sigma, M)$ and
the Lipschitz character of $\Omega$.
By covering $\partial\Omega$ with balls centered on $\partial\Omega$
with radius $ c\e$, we may deduce from (\ref{8.1-2}) that
\begin{equation}\label{8.1-3}
\int_{\partial\Omega} |\nabla u_\e|^2\, d\sigma
\le C \int_{\partial\Omega} |g|^2\, d\sigma
+\frac{C}{\e} \int_{\Omega_{\e}} |\nabla u_\e|^2\, dx.
\end{equation}
Since $u_\e \perp \mathcal{R}$ in $L^2(\Omega; \rd)$,
it follows from Theorem \ref{large-scale-Rellich-theorem} that
$$
\int_{\partial\Omega} |\nabla u_\e|^2 \, d\sigma \le C \int_{\partial\Omega} |g|^2\, d\sigma.
$$
Thus we have proved the Rellich estimate (\ref{Rellich-8}) and Theorem \ref{Rellich-theorem}.
\end{proof}

\bibliography{Geng-Shen-Song-2016.bbl}

\medskip

\begin{flushleft}
Jun Geng,
School of Mathematics and Statistics,
Lanzhou University,
Lanzhou, P.R. China
\quad

E-mail:gengjun@lzu.edu.cn
\end{flushleft}

\begin{flushleft}
Zhongwei Shen,
 Department of Mathematics,
University of Kentucky,
Lexington, Kentucky 40506,
USA. 

E-mail: zshen2@uky.edu
\end{flushleft}

\begin{flushleft}
Liang Song, 
Department of Mathematics,
Sun Yat-sen University,
Guanzhou, P.R. China

E-mail: songl@mail.sysu.edu.cn

\end{flushleft}
\medskip

\noindent \today

\end{document}